\input psfig
\def\endofproof{ \hfill $\#$ }
\def\wlog{ without loss of generality }

\def\hf{Hubbard Forest }
\def\Hf{Hubbard Forest}
\vskip 1.5in
\bigskip
\bigskip
\centerline{\bf Coexistence of Critical Orbit Types}
\centerline{\bf in Sub-Hyperbolic Polynomial Maps}
\bigskip
\bigskip
\centerline{Alfredo Poirier}
\centerline{Institute for Mathematical Sciences}
\centerline{SUNY at StonyBrook}
\centerline{StonyBrook, NY 11794}
\bigskip
\bigskip
{\bf Abstract.} 
{\sl We establish necessary and sufficient conditions for the realization of mapping schemata 
as post-critically finite polynomials, 
or more generally, 
as post-critically finite polynomial maps 
from a finite union of copies of the complex numbers {\bf C} to itself 
which have degree two or more in each copy. 
As a consequence of these results we prove a transitivity relation 
between hyperbolic components in parameter space 
which was conjectured by Milnor.} 
\vskip 0.5in
\centerline{\bf 1. Introduction.}
\bigskip
\bigskip
If  $f$  is a proper holomorphic map of the complex numbers  ${\bf C}$, 
or more generally from a finite union of copies of  ${\bf C}$ to itself, 
then the proper homotopy class of  $f$  can be described by a very simple
combinatorial structure which will be called the ``ambient mapping
schema'' ${\cal M}$ 
(see the discussion below, Appendix B, and compare [M]). 
If  $f$  is post-critically finite, 
then the restriction of  $f$  to the union of critical orbits is described by another mapping schema  $S$ 
(compare Appendix B). 
There is a natural ``projection map'' from  $S$  to  ${\cal M}$,
which sends each post-critical point to the copy of ${\bf C}$ containing it.
The object of this note is to characterize exactly which projection maps
$S \to {\cal M}$  can arise in this way.
\bigskip
{\bf 1.1 Definition: Mapping Schemata. } 
By a {\it mapping schema} $S=(|S|,F,w)$ (or a schema in short) we mean: 

(1) a finite set $|S|$ of points, together with 

(2) a function $F$ from $|S|$ to itself, and also 

(3) a ``weight function" $w$ which assigns an integer $w(v) \ge 0$ called the {\it critical weight} 
to each $v \in |S|$. 
\medskip\noindent
Equivalently, 
such a mapping schema can be represented by a finite graph with one vertex for each $v \in |S|$, 
and with exactly one directed edge $e_v$ leading out from each vertex 
$v$ to a vertex $F(v)$. 
By definition the degree associated with the edge $e_v$ 
(or with the vertex $v$) is the integer $d(v)=w(v)+1 \ge 1$. 
The {\it weight $w(S)$ of a schema} $S$ is by definition the number 
$w(S)=\sum_{v \in |S|} w(v)$, 
the {\it degree of $S$} is then defined as $deg(S)=w(S)+1$. 
\eject
\medskip
Such a mapping schema is {\it reduced} if every vertex is critical. 
Suppose that we start with a mapping schema $S$ which satisfies the following very mild condition: 
{\it 
Every cycle in S contains at least one critical vertex.} 
Then there is an 
{\it associated reduced mapping schema} $\bar S$ 
which is obtained from $S$ simply by discarding all vertices of weight zero and shrinking every edge of degree one to a point. 
Note that $S$ and $\bar S$ have the same total weight. 
\bigskip 
{\bf 1.2 Definition: Ambient Schemata. } 
Let ${\cal M}=(|{\cal M}|,F_{\cal M},w_{\cal M})$ be a schema such that every cycle of ${\cal M}$ contains at least a critical vertex. 
Form the disjoint union ${\cal M}\times {\bf C}$ of $n$ copies of the complex numbers ${\bf C}$, 
where $n$ is the number of vertices of ${\cal M}$. 
In other words, 
replace each vertex $u \in {\cal M}$ by a copy of ${\bf C}$. 
Let ${\cal P}^{\cal M}$ be the space consisting of all maps $f$ from $|{\cal M}|\times {\bf C}$ to itself 
such that the restriction of $f$ to each component $u \times {\bf C}$ is a 
monic centered polynomial of degree $d_{\cal M}(u)=w_{\cal M}(u)+1$, 
taking values in $F_{\cal M}(u) \times {\bf C}$. 
By definition then ${\cal M}$ is the {\it ambient schema} of $f$. 
Note that whenever $f$ is a proper holomorphic map from a finite union of copies of  ${\bf C}$, to itself, 
then after component-wise affine change of coordinates, 
necessarily $f \in {\cal P}^{\cal M}$ for some schema ${\cal M}$. 
\smallskip
{\bf Remark.} 
Note that in the definition above, 
whenever $\overline{\cal M}$ is the associated reduced schema of ${\cal M}$, 
there is a canonical bijective correspondence between maps 
$f \in {\cal P}^{\cal M}$ and maps $\bar f \in {\cal P}^{\overline{\cal M}}$ 
defined as follows. 
For $f \in {\cal P}^{\cal M}$ we obtain $\bar f$ by simply discarding all the identity maps. 
Conversely from $\bar f$ we can obtain a unique $f$ by interpolating the identity in the adequate places. 
\bigskip
{\bf 1.3 Definition: Post-critical Schemata.} 
Let $f \in {\cal P}^{\cal M}$ be post-critically finite 
(i.e, every critical point of $f$ eventually maps to a periodic cycle). 
Denote by $\Omega(f)$ the set of critical points of $f$. 
Let $|S|$ be any finite invariant set containing this critical set $\Omega(f)$. 
This $f$ and $|S|$ define a schema $S=(|S|,f,w_f)$, 
where $deg_f(v)=w_f(v)+1$ is the local degree of $f$ at $v \in |S|$. 
By definition this schema is called a {\it post-critical schema}. 
(Note that this terminology may not be standard. 
In fact, we will use the word post-critical even if we allow vertices which do not belong to the orbit of this critical set. 
We reserve the word {\it ``minimal"}
for that schema generated by the critical set.) 
Again ${\cal M}$ can be thought as the ambient schema of $S$ in the sense described in the next definition. 
For the relations between these several definitions, 
the reader is referred to Appendix B. 
\bigskip
{\bf 1.4 Definition: Projection Maps between Schemata.} 
Let $S=(|S|,F_S,w_S)$ and ${\cal M}=(|{\cal M}|,F_{\cal M},w_{\cal M})$ be schemata. 
By a {\it projection map} $\phi$ from $S$ to ${\cal M}$ will be meant a map 
$\phi:|S| \to |{\cal M}|$ which satisfies the following conditions.
\smallskip
a) {\it $\phi$ semiconjugates $F_S$ to $F_{\cal M}$;} 
in other words the diagram below commutes
$$
\matrix{
&\hbox{ }&|S|&{F_S \atop \longrightarrow} &|S| &\hbox{ } \cr
&\phi &\Big\downarrow & &\Big\downarrow  &\phi \cr  
&\hbox{ }&|{\cal M}|&{F_{\cal M} \atop \longrightarrow} &|{\cal M}| &\hbox{ } \cr
}
$$
\smallskip
b) {\it $\phi$ preserves weight,} in the sense that for each $u \in |{\cal M}|$ 
$$ 
w_{\cal M}(u)=\sum_{\phi(v)=u} w(v).
$$
\smallskip
If this is the case and every cycle of ${\cal M}$ contains a critical vertex, 
${\cal M}$ will be called 
{\it the ambient schema of $S$}. 
(This terminology is an abuse of language; 
${\cal M}$ should be really called the ambient schema of $S$ under $\phi$.) 
For $u \in |{\cal M}|$ we set $W(u)=\{v \in |S|: \phi(v)=u \}$ and call it the {\it fibre at $u$}. 
\bigskip 
{\bf 1.5 Remark.} 
In practice this last definition is usually given in a constructive way.\break 
In general we will start with a schema $S=(|S|,F,w)$ 
(which we may think as a map which is component-wise polynomial) 
and a function 
$F_{\cal M}:|{\cal M}| \to |{\cal M}|$ from a finite set to itself 
(which we may think as the pattern in which each copy of the complex numbers maps to another). 
Whenever a ``projection map" 
$\phi:|S| \to |{\cal M}|$ is given, 
is onto and semiconjugates $F_S$ to $F_{\cal M}$, 
we can canonically define ``an ambient schema" ${\cal M}=(|{\cal M}|,F_{\cal M},w_{\cal M})$ by setting 
$$
w_{\cal M}(u)=\sum_{\phi(v)=u} w_S(v).
$$
If every cycle in $|{\cal M}|$ contains a critical point 
(i.e, $w_{\cal M}(u) \ge 1$ for some $u$ in the cycle), 
then by definition {\it $S$ projects to ${\cal M}$}. 
In fact, in practice we might be given the critical orbits of a collection of polynomial maps $f$; i.e, a schema $S$. 
If furthermore we know to which copy of the complex numbers each point belongs, 
we can easily reconstruct the unique schema ${\cal M}$ for which 
$f \in {\cal P}^{\cal M}$. 
\bigskip
Of course, if ${\cal M}$ is a reduced schema and  
$f \in {\cal P}^{\cal M}$ is post-critically finite; 
then every associated post-critical schema $S_f$
projects to ${\cal M}$ by definition. 
The purpose of this note is to 
characterize those mapping schemata which can arise in this way. 
In other words, which schemata can be realized as a 
post-critical schema of a sub-hyperbolic polynomial in the space ${\cal P}^{\cal M}$. 
Note that for this, 
we only have to study that case in which 
the ambient schema consists of a single connected component. 
To state these conditions (for realization) we will need some preliminary definitions. 
\bigskip
{\bf 1.6 Definition.} 
Fix an $n \ge 2$, and consider the polynomial $P_n(z)=z^n$. 
Let $N(n,k)$ be by definition the number of disjoint periodic orbits of period (strictly) $k$ under iteration of this polynomial $P_n$. 
The following result is folklore. 
\bigskip
{\bf Theorem.} 
{\it Let $P$ be a post-critically finite polynomial of degree $n \ge 2$. 
Then the number of disjoint periodic orbits of period $k$ is $N(n,k)$.} 
\medskip
{\bf Proof.} 
The proof is left to the reader, 
who has to recall that the equation 
$P^{\circ k}(z)=z$ has only simple solutions in this post-critically finite case. 
\endofproof
\bigskip
\bigskip
{\bf Remark.} 
For $n \ge 2$ note that $n-1 \le N(n,k)$; with equality holding if and only if 
$n=k=2$. 
\bigskip
If the ambient schema ${\cal M}$ consists of a single connected component, 
we denote by $n({\cal M})$ the product of the degrees of all elements which belong to the unique cycle contained in ${\cal M}$. 
By hypothesis this {\it inner degree of ${\cal M}$} is always bigger than 1. 
\bigskip
{\bf 1.7 Definition.} 
Suppose $\phi$ is a projection from $S$ to the connected ambient\break 
schema ${\cal M}$. 
Now let $v \in |S|$ be periodic; 
then $F=F_S$ induces a periodic map in the fibre  
$W(\phi(v)) \cap \{v,F(v),F^{\circ 2}(v),\dots\}$. 
By definition this period $k$ will be called {\it the return period of $v$}. 
In practice, the return period can be computed by dividing the `actual' period of $v$ by the number of periodic points in the connected ambient schema ${\cal M}$. 
\bigskip
{\bf Theorem A.} 
{\it Let ${\cal M}$ be a connected ambient schema with a critical cycle.\break 
Let $f \in {\cal P}^{\cal M}$ be post-critically finite, 
and $|S|$ a finite invariant set containing all critical points of $f$. 
Then the canonical projection $\phi$ from the post-critical schema $S$ to the ambient schema ${\cal M}$ is admissible the sense that 
\smallskip
a) For every $u \in |{\cal M}|$ and $v \in W(F_{\cal M}(u))$ we have 
$$
\sum_{\{v' \in W(u), F_S(v')=v\}} d_S(v') \le d_{\cal M}(u).
$$
\smallskip
b) The number of disjoint periodic orbits in $S$ of return period $k$ 
is less than or equal to $N(n({\cal M}),k)$.} 
\bigskip
{\bf Proof.} 
Let $f_u=f|_{u \times {\bf C}}$. 
This $f_u$ is a polynomial of degree $d_{\cal M}(u)$ which maps $u \times {\bf C}$ to $F_{\cal M}(u) \times {\bf C}$. 
Therefore every point $v \in W({F_{\cal M}(u)})$ has at most $d_{\cal M}(u)$ inverses counting multiplicity. 
This establishes property a). 
Furthermore, 
if $u_0 \mapsto u_1 \mapsto \dots \mapsto u_{r}=u_0$ is the cycle in ${\cal M}$, 
then the post-critically finite polynomial $P=f_{u_{r-1}}\circ \dots \circ f_{u_0}$ has degree $n({\cal M})$. 
The result then follows from Theorem 1.6. 
\endofproof
\bigskip
{\bf 1.8 Definition: Admissible Projections.} 
The previous result motivates the following definition. 
Let $\phi$ be a projection from $S$ to ${\cal M}$. 
Suppose first the ambient schema ${\cal M}$ is connected and its only cycle contains at least a critical vertex.  
We say that $\phi$ is {\it admissible} 
(or in short that $S$ is admissible on ${\cal M}$) if 
{\it 
\smallskip
a) For every $u \in |{\cal M}|$ and $v \in W(F_{\cal M}(u))$ we have 
$$
\sum_{\{v' \in W(u), F_S(v')=v\}} d_S(v') \le d_{\cal M}(u).
$$
\smallskip
b) The number of disjoint periodic orbits in $S$ of return period $k$ 
is less than or equal to $N(n({\cal M}),k)$.}
\bigskip
More generally, an ambient schema ${\cal M}$ can always be described 
as the disjoint union of connected schemata. 
Therefore $\phi$ can be decomposed into projections to different connected 
ambient schemata. 
In this case we say that $\phi$ is {\it admissible} 
if all such projections are admissible in the sense described above. 
\bigskip
\bigskip
\eject
{\bf Theorem B. } 
{\it Conversely, if $\phi$ defines an admissible projection from 
a mapping schema $S$ to an ambient schema ${\cal M}$,   
then there is a post-critically finite map $f \in {\cal P}^{\cal M}$ and a finite invariant set containing the critical points of $f$ whose associated post-critical schema is isomorphic to $S$.} 
\bigskip
The proof of this Theorem will be given starting in Section 3. 
Clearly we shall restrict ourselves only to connected ambient schemata. 
\bigskip
\bigskip
{\bf 1.9 Caution.} 
That condition in the bound in the number of periodic cycles of given period 
is important even in the case 
where there are no `superfluous cycles' in $S$; that is, 
in the case where every connected component of $S$ contains a critical point. 
For example, consider the following schema: 
\medskip
\centerline{\psfig{figure=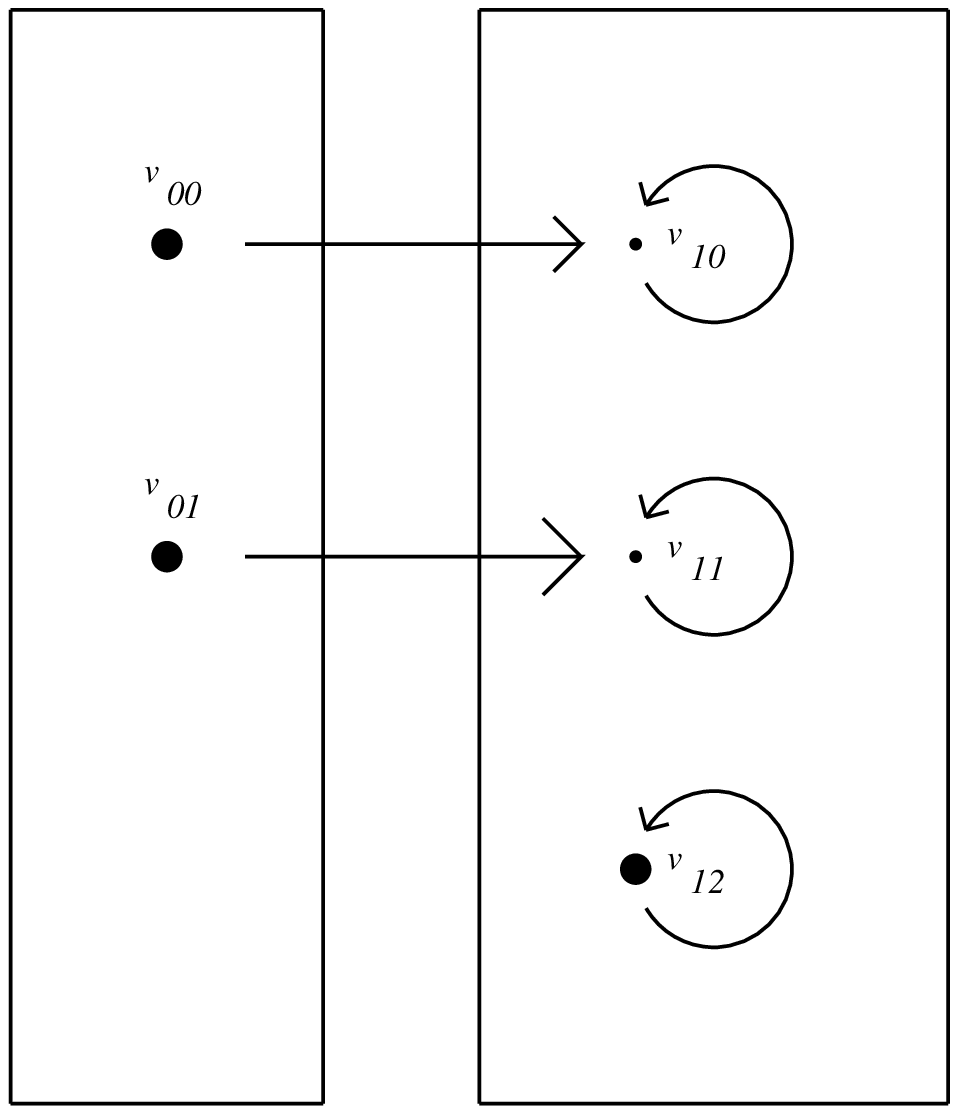,height=3.5in}}
\medskip
\noindent 
Define $F(v_{0i})=v_{1i}$ and $F(v_{1i})=v_{1i}$. 
Let $v_{00},v_{01},v_{12}$ be the unique (simple) critical points.  
Take as projection map $\phi(v_{ij})=i$; 
and note that this projection is not admissible. 
In fact, if this schema were to be realized, 
we will have a polynomial map of degree $2$ with 3 fixed points, which is impossible. 
\bigskip 
However, in case the ambient schema ${\cal M}$ is cyclic, 
it is easy to see that in the absence of ``superfluous cycles" condition b) is implied by condition a). 
In fact, in this case 
every component of $S$ contributes at least with one to the degree of the inner part of ${\cal M}$. 
Therefore the allowed number of cycles of a given period is also increased (compare Theorem C). 
On the other hand, if ${\cal M}$ is not cyclic 
(or if ``superfluous components" are allowed) 
then a 
component of $S$ does not necessarily contribute to the degree of the inner part of ${\cal M}$. 
Thus we can be lead to an ``inflationary" process in the number of cycles of a given period if condition b) is relaxed (compare $\S$1.9 above). 
\bigskip
{\bf 1.10 Definition.} 
A {\it semi-reduced schema} $S=(|S|,F,w)$ is a schema satisfying the following conditions. 
\smallskip
a) {\it Every connected component of $S$ contains a critical point;} and 
\smallskip
b) {\it Every end of $S$ is critical.} 
\bigskip
{\bf Remark.} 
Given $f \in {\cal P}^{\cal M}$ post-critically finite with critical set $\Omega(f)$, 
we have a canonically defined semi-reduced schema 
$S_f=(|S_f|={\cal O}(\Omega(f)),P,w_f)$. 
Here ${\cal O}(\Omega(f))$ is the orbit of the critical set. 
This mapping schema can be view as the 
``minimal post-critical schema of $f$". 
\bigskip
\bigskip
{\bf Theorem C.} 
{\it Suppose the semi-reduced schema $S$ projects to a cyclic ambient schema ${\cal M}$.
Then there is a post-critically finite map $f \in {\cal P}^{\cal M}$ 
with associated minimal post-critical schema isomorphic to $S$ 
if and only if 
for every $u \in |{\cal M}|$ and $v \in W(F_{\cal M}(u))$ we have 
$$
\sum_{\{v' \in W(u), F_S(v')=v\}} d_S(v') \le d_{\cal M}(u).
$$
}
\smallskip
{\bf Proof.} 
Suppose $S$ consists of $n-1$ components. 
As each of these components contains a critical point it follows easily that
$n \le deg(S)$. 
By Lemma 1.11 below, we have $deg(S) \le n({\cal M})$; 
and therefore $n \le deg(S) \le n({\cal M})$. 
Now the result follows from Theorem B and the fact that 
$n-1 \le N(n,k) \le N(n({\cal M}),k)$ for all $k$. 
\endofproof
\bigskip
{\bf 1.11 Lemma.} 
{\it Suppose the semi-reduced schema $S$ projects to a cyclic ambient schema ${\cal M}$.
Then $deg(S) \le n({\cal M})$.} 
\medskip
{\bf Proof.} 
By definition 
$$
deg(S)=1+\sum_{v \in \phi^{-1}(|{\cal M}|)} w_S(v)=1+\sum_{u \in |{\cal M}|} \sum_{v \in \phi^{-1}(u)} w_S(v), 
$$
and 
$$
n({\cal M})=\prod_{u \in |{\cal M}|} (1+w_{\cal M}(u))=\prod_{u \in |{\cal M}|} (1+\sum_{v \in \phi^{-1}(u)} w_S(v)).
$$
(Compare condition b) in Definition 1.4.) 
The result follows then by elementary properties of non negative real numbers. 
\endofproof
\bigskip
Theorem C is true in particular for post-critically finite polynomials: 
\bigskip
{\bf Theorem D.} 
{\it A semi-reduced schema $S$ of degree $n=deg(S) \ge 2$ can be realized by a polynomial of degree $n$ if and only if  
for every $v \in |S|$ we have $\sum_{v' \in |S|:F(v')=v} d_S(v') \le n$.} 
\bigskip
{\bf Proof.} 
As $S$ consists of at most $n-1$ components, 
the result follows from Theorem B and the fact that 
$n-1 \le N(n,k)$ for all $k$. 
\endofproof
\bigskip
\bigskip
As an application of Theorem B, 
we prove in Section $2$ a conjecture of Milnor\break 
(compare [M, remark 2.11]). 
For two ambient schemata ${\cal M}$ and ${\cal M}'$ of the same weight 
we write ${\cal M} \succ {\cal M}'$ if and only if 
the connectedness locus ${\cal C}({\cal M})$ 
contains a hyperbolic component with reduced schema isomorphic to ${\cal M}'$. 
We prove that this relation between ambient schemata is transitive 
(Corollary 2.7). 
Jeremy Kahn has pointed out that this result can also be proved using Levy-Bernstein Theorem which states that: 
{\it a marked\break 
topological polynomial which has the property that every critical point eventually falls into a critical cycle is 
Thurston equivalent to a unique post-critically finite polynomial.} 
\medskip
We will prove Theorem B starting in Section 3. 
In essence, our proof  can be taken as constructive. 
We will ``construct" in each case a \hf 
(of the appropriate degree) 
in which a suitable set of ``marked vertices" will realize the given schema. 
Our strategy will be developed in Section 3, 
where we also include some useful remarks. 
Then we continue in Sections 4 and 5 by realizing the easist possible cases. 
In Section 6 we show how to modify a \hf 
which realizes a given schema, 
in order to realize a bigger one on which the former is properly contained. 
The case in which a component of the schema contains a periodic cycle of `return period $2$' 
represents the hardest isolated case and requires special treatment 
(compare Example 3.3 and Lemma 7.7). 
This can be explained by the unfortunate fact that $N(2,2)=1$. 
In Section $8$ we proof Theorem B for the special case in which the ambient schema is cyclic. 
After this, the proof of Theorem B in the general case is routine and is given in Section 9. 
We also include Appendix A where we state the main definitions and results related to \Hf s; 
and Appendix B, where we state the relations among the several possible 
interpretations of schemata in Complex Dynamics.  
\bigskip
{\bf Acknowledgement.} 
We will like to thank John Milnor for helpful conversations. 
Also, we want to thank the Geometry Center, University of Minnesota
and Universidad Cat\'olica del Per\'u 
for their material support. 
\bigskip
\bigskip
\centerline{\bf 2. An Application to Parameter Space.} 
\bigskip
\bigskip
{\bf 2.1 Definition.} 
The {\it connectedness locus} ${\cal C}({\cal M})\subset {\cal P}^{\cal M}$ 
is by definition the set of maps $f \in {\cal P}^{\cal M}$ for which 
the orbit ${\cal O}(\Omega(f))$ of the critical set remains bounded under iteration. 
Equivalently $f \in {\cal C}({\cal M})$
if and only if 
the restriction of the filled Julia set $K(f)$ 
to each copy $u \times {\bf C}$ of the complex numbers is connected; 
in other words, 
if the set of points of bounded orbit is connected in each copy of ${\bf C}$. 
\medskip
Every hyperbolic component ${\cal H} \subset {\cal C}({\cal M})$ 
contains a unique center point $f$ which is post-critically finite. 
In this case, 
all critical points of $f$ belong to the Fatou set $F(f)$. 
Thus, this $f$ has a canonical 
``{\it minimal post-critical  mapping schema}", as defined by\break 
$S_f=(|S_f|={\cal O}(\Omega(f)), f, w_f)$. 
It follows that $S_f$ is semi-reduced 
(compare Definition 1.10). 
Also, because $f$ is post-critically finite and hyperbolic, 
every connected component of $S_f$ contains a critical cycle. 
\bigskip
{\bf 2.2 Definition.} 
By the {\it type of a hyperbolic component ${\cal H}$} will be meant the reduced schema 
$\overline{S_f}$ associated to the minimal post-critical schema $S_f$ of the center of ${\cal H}$. 
It is proved in [M] Theorem 5.1, that two hyperbolic components 
(regardless of the parameter space on which they are defined) 
are biholomorphically equivalent if 
they have isomorphic mapping schemata. 
\medskip
{\bf Remark.} 
Note that by definition the mapping schema $S_f$ projects to ${\cal M}$, 
while in general there is no relation between $\overline{S_f}$ and ${\cal M}$.  
\bigskip
{\bf 2.3 Definition.} 
If the connectedness locus ${\cal C}({\cal M})$ contains a hyperbolic component of 
type ${\cal M}'$ we write ${\cal M} \succ {\cal M}'$. 
\bigskip
{\bf 2.4 Caution:} 
It may happen that ${\cal M} \succ {\cal M}'$ and ${\cal M}' \succ {\cal M}$, 
even though ${\cal M}$\break 
is not isomorphic to ${\cal M}'$. 
For example this is true for ${\cal M}=(\{2\}1,\{1\}1,1)$ and\break 
${\cal M}'=(\{2\}1,1,\{1\}1)$ with diagrams
\smallskip
\centerline{\psfig{figure=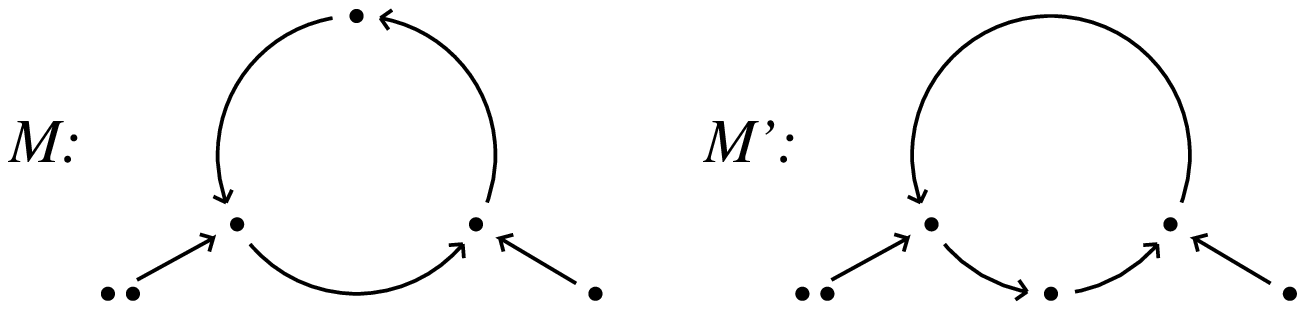,height=1.5in}}
\medskip
Figure 1 shows a ``Hubbard Forest" which realizes ${\cal M} \succ {\cal M}'$, 
and similarly we can show that ${\cal M}' \succ {\cal M}$. 
\bigskip 
\medskip
\centerline{\psfig{figure=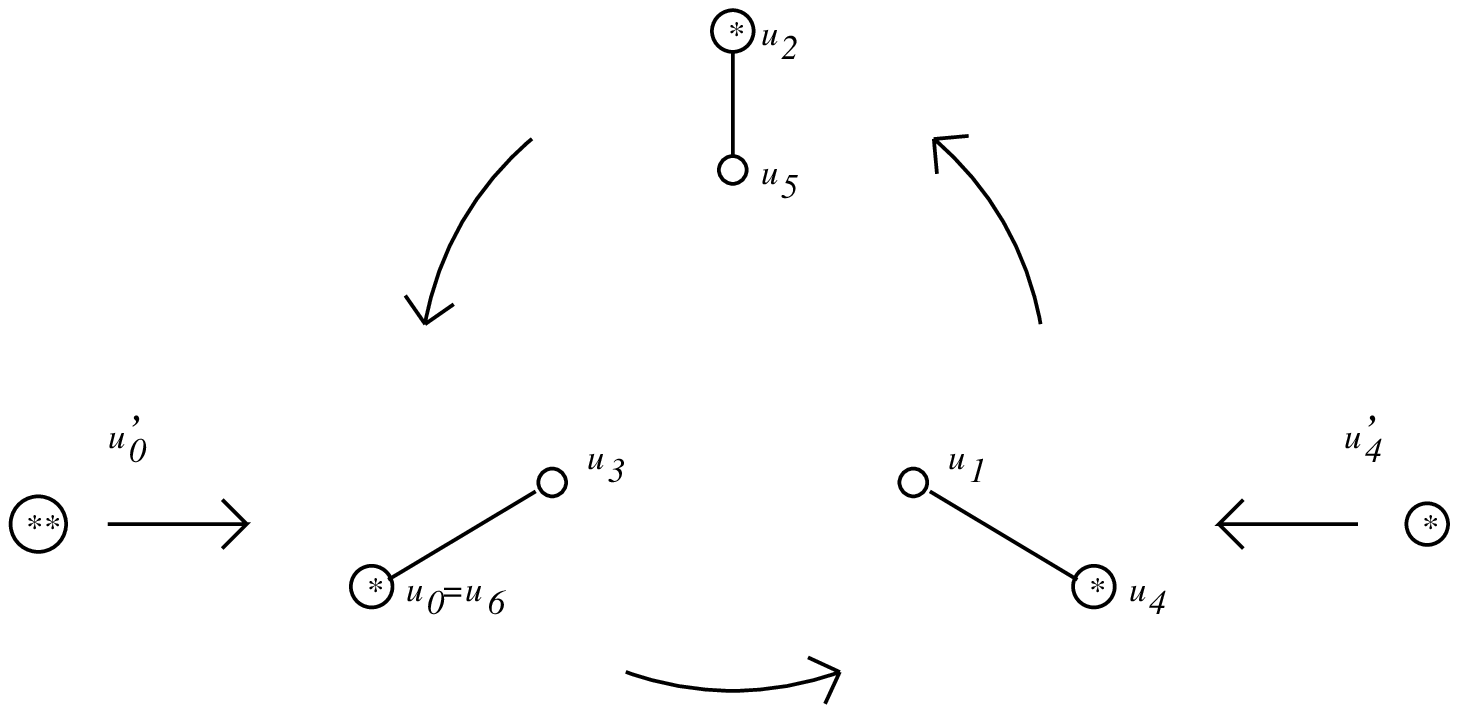,height=2.5in}}
\noindent
{\bf\sl Figure 1. A Hubbard Forest which shows that ${\cal M} \succ {\cal M}'$. 
This Hubbard Forest represets the center of a hyperbolic component in ${\cal C}({\cal M})$ with reduced schema isomorphic to ${\cal M}'$.\break 
Here 
$u_0 \mapsto u_1 \mapsto u_2 \mapsto u_3 \mapsto u_4 \mapsto u_5 \mapsto u_0$ 
and $u'_0 \mapsto u_0$, $u'_4 \mapsto u_4$.} 
\bigskip
\eject
{\bf 2.5 Theorem.}  
{\it The connectedness locus ${\cal C}({\cal M})$ 
contains a hyperbolic component with\break 
reduced schema isomorphic to ${\cal M}'$ (i.e, ${\cal M} \succ {\cal M}'$)  
if and only if 
there is a semi-reduced schema $S$ satisfying the following two conditions 
\smallskip
a) The associated reduced schema $\bar S$ is well defined and isomorphic to ${\cal M}'$. 
\smallskip
b) There is an admissible projection from $S$ to ${\cal M}$.} 
\medskip
{\bf Proof.} 
This follows easily from the definition and Theorem B. 
In fact, 
if ${\cal M} \succ {\cal M}'$ then by definition there is a hyperbolic component 
${\cal H}$ of type ${\cal M}'$ with center $f$.  
Then we can take $S$ as the minimal post-critical schema of $f$. 
Conversely, 
if there is such a semi-reduced schema $S$, then Theorem B guarantees that there is an associated post-critically finite $f \in {\cal P}^{\cal M}$. 
Clearly this $f$ is hyperbolic and the connected hyperbolic component with center $f$ is of type ${\cal M}'$. 
\endofproof
\medskip
{\bf 2.6 Remark.} 
Our main goal of this Section is to prove that this relation between 
ambient schemata is transitive. 
For this it is important to recall that whenever 
the schema $S$ has associated reduced schema ${\cal M}$, 
there is a canonical isomorphism between the spaces ${\cal P}^{\cal S}$ and ${\cal P}^{\cal M}$. 
Furthermore, if $f_{\cal M}$ is hyperbolic and postcritically finite, 
then the corresponding $f_{S}$ has the same properties and 
the associated reduced mapping schema of both functions are the same. 
\bigskip
{\bf 2.7 Corollary.} 
{\it This relation $\succ$ between reduced schemata is transitive.} 
\medskip
{\bf Proof.} 
Suppose ${\cal M}_0 \succ {\cal M}_1$ and ${\cal M}_1 \succ {\cal M}_2$. 
We may assume that ${\cal M}_0$ is connected. 
Then there are hyperbolic polynomials maps $f_1 \in {\cal P}^{{\cal M}_0}$, 
$f_2 \in {\cal P}^{{\cal M}_1}$ whose semireduced schemata $S_i$ ($i=1,2$) 
have associated reduced schemata isomorphic to ${\cal M}_i$. 
In this way, by the previous remark, there is a polynomial map 
$\tilde f_2 \in {\cal P}^{S_1}$ whose reduced schema ${\cal S}$ is isomorphic to ${\cal M}_2$. 
To complete the proof we only have to construct an admissible projection from 
${\cal S}$ to ${\cal M}_0$. 
But this is easily done by composing the admissible projections $\psi_1$ from $S_1$ to ${\cal M}_0$ with $\psi_2$ from ${\cal S}$ to $S_1$. 
Clearly this composition semiconjugates $F_S$ to $F_{{\cal M}_0}$ and preserves weights. We still need to prove the admissibility conditions. 
\smallskip
Let $n$ denote the sum of the weights of the periodic points in ${\cal M}_0$. 
As every periodic orbit in ${\cal S}$ contains at least a critical point
(which furthermore can only project to a periodic point in ${\cal M}$), 
it follows that  
the number of periodic orbits in ${\cal S}$ is less or equal to $n$. 
In particular the number of disjoint periodic orbits of return period $k$ 
is always at most $n$, which in turn is bounded by $N(n({\cal M}),k)$. 
(Compare Lemma 1.11 and the proof of Theorem C.) 
\smallskip
Now let $u \in |{\cal M}_0|$ and $v \in (\psi_1\circ \psi_2)^{-1}(F_{{\cal M}_0}(u))$. 
To simplify notation we define\break 
the following sets. 
Let ${\cal X}=\{{v' \in (\psi_1 \circ \psi_2)^{-1}(u): F_{\cal S}(v')=v}\}$
and for any\break 
$w \in {\cal W}=\{w' \in \psi_1^{-1}(u): F_{S_1}(w')=\psi_2(v)\}$ 
define 
${\cal Y}_w=\{v' \in \psi_2^{-1}(w): F_{\cal S}(v')=v\}$. 
Then clearly ${\cal X}=\cup_{w \in {\cal W}} {\cal Y}_w$. 
It follows that 
$$
\sum_{v' \in X} d_{\cal S}(v')=\sum_{w \in {\cal W}} \sum_{v' \in {\cal Y}_w} 
d_{\cal S}(v') \le \sum_{w \in {\cal W}} d_{S_1}(w) \le d_{{\cal M}_0}(u); 
$$
where the last two inequalities follow from the admissibility 
(in the sense of Theorem A) of $\psi_1$ and $\psi_2$.
The result now is a consequence of Theorem 2.5. 
\endofproof
\bigskip
\bigskip
\vfill\eject
 
\centerline{\bf 3. Strategy.} 
\bigskip
\bigskip
We will proof Theorem B by showing the existence of 
a Hubbard Forest ${\bf H^*}(S)$ on which the schema $S$ can be realized. 
(For the convenience of the reader we have included in the appendix 
the notation and results from [P3] where this concept was introduced.) 
In other words, 
let $S$ be admissible on ${\cal M}$, that is, 
suppose there is an admissible projection map $\phi:|S| \to |{\cal M}|$. 
To realize $S$ by a map $f \in {\cal P}^{\cal M}$ it is enough to 
find a \hf 
${\bf H^*}$ with set of vertices $V$ 
which realizes a post-critically finite map in ${\cal P}^{\cal M}$ 
with certain obvious properties
(in particular there is a canonical projection from $V$ to ${\cal M}$).  
This will happen if the mapping schema $S$ and the forest ${\bf H^*}$ are related in the following way: 
\medskip
{\it 
There is an embedding $\alpha$ from the set $|S|$ to the set $V$ of vertices of the forest which 
commutes with the dynamics and preserves weight.  
We require that both, a vertex $v$ in $|S|$ and its embedded image $\alpha(v) \in V$ in the forest, 
project to the same $u \in |{\cal M}|$.}  
\medskip
A point of the form $\alpha(v)$ is a {\it marked vertex (respect to $S$)}. 
To simplify notation, 
we will not distinguish between a point in $|S|$ and its image in the forest
unless strictly necessary. 
When there is no possible confusion we will simply say that $v$ is a {\it marked vertex.}  
\bigskip 
We have thus reduced the proof of Theorem B to the construction of a purely combinatorial object. 
In most of the discussion that follows, 
it is convenient to think of $S$ as being semi-reduced 
(compare Definition 1.10). 
However, note that unless clearly specified, 
all statements and definitions apply as well to more general admissible schemata. 
\bigskip
The simplest mapping schema we can think of, 
is that consisting of a critical cycle. 
This schema $S$ necessarily projects to a cyclic ambient schema and is admissible. 
In fact, 
that the ambient schema is cyclic follows from the fact that $S$ is cyclic and the projection is onto. 
That the projection is admissible can be easily seen as follows. 
We start by assigning a ``degree" equal to one to every $u \in |{\cal M}|$. 
Then each vertex $v \in |S|$ contributes with an additional $w_S(v)=d_S(v)-1$ to the degree of $\phi(v) \in |{\cal M}|$; 
while it contributes with $d_S(v)$ to the sum 
$\sum_{\{v' \in W(\phi(v)):F_S(v')=F_S(v)\}} d_S(v')$. 
As in this case every vertex $v \in |S|$ has exactly one preimage, 
the claim follows. 
We will realize cyclic schemata 
(and more general disjoint union of cyclic schemata) in Section 4. 
\bigskip
The first technical difficulties can only appear if a given vertex in $S$ has two distinct preimages in the same `fibre' $W(u)$. 
We start studying this case with the following lemma. 
\bigskip
\eject
{\bf 3.1 Lemma.} 
{\it Suppose $S$ is admissible on ${\cal M}$, 
and let 
$v_1,v_2 \in W(u)$ be different.\break 
If $F(v_1)=F(v_2)$ 
then there exists $v \in W(u)-\{v_1,v_2\}$ such that 
$d_S(v) > 1$. 
Furthermore, this critical vertex $v$ can be chosen so that $F(v) \ne F(v_1)$.} 
\medskip
{\bf Proof.} 
Let $v_1,\dots v_n$ (with $n>1$) 
be the maximal collection of vertices such that $F(v_i)=F(v_1)$. 
Then 
$$
\sum_{i=1}^{n} d_S(v_i) = \sum_{v \in W(u):F(v)=F(v_1)} d_S(v) 
\le d_{\cal M}(u)=1+\sum_{v \in W(u)} (d_S(v)-1). 
$$
In this way there is a vertex $v \in W(u)$ other than $v_1, \dots v_n$, 
for which $d_S(v)-1>0$. 
\endofproof
\bigskip
This last result is not surprising if we adopt another point of view. 
In fact, if $S$ could be realized by a Hubbard Forest ${\bf H^*}$, 
then $v_1$ and $v_2$ should belong to the same connected component of ${\bf H^*}$. 
As they map to the same point, 
it follows that there is always an interior critical point in the path $[v_1,v_2]$ joining these vertices. 
We have proved the following. 
\bigskip
{\bf 3.2 Lemma.} 
{\it Suppose the Hubbard Forest ${\bf H^*}$ realizes the projection from $S$ to ${\cal M}$. 
Let $v_1,v_2 \in W(u)$ be such that $v_1 \ne v_2$  
and denote by $H_u$ the connected component of ${\bf H^*}$ to which $v_1,v_2$ belong. 
If $F_S(v_1)=F_S(v_2)$ 
then there exists $v \in W(u)-\{v_1,v_2\}$ such that 
$d_S(v) > 1$ and 
$v_1,v_2$ belong to different components of $H_u-\{v\}$. 
Furthermore, this critical vertex $v$ can be chosen so that $F(v) \ne F(v_1)$.} 
\endofproof
\bigskip
{\bf 3.3 Example: Critical points mapping into the same cycle.} 
Let $S$ be any schema which consist of a cycle ${\cal C}$ 
and two critical points $w_0 \ne w_1$ which not belong to the cycle but 
such that $F(w_0) \ne F(w_1)$ both belong to the cycle. 
In particular ${\cal C}$ should have period at least $2$. 
In other words $S$ can be described as follows: 
This $S$ consists of a cycle 
${\cal C}: v_0 \mapsto v_1 \mapsto \dots \mapsto v_m=v_0$ 
of period at least two. 
This cycle ${\cal C}$ may contain critical vertices or not. 
Also there two vertices $\omega_0,\omega_1$ outside this cycle; 
that is, $S={\cal C} + \{\omega_0,\omega_1\}$. 
We require both $\omega_0$ and $\omega_1$ to be critical and such that 
$F(\omega_0)=v_1 \ne F(\omega_1)=v_j$ for some $v_1$ and $v_j$ different elements in the cycle. 
(Compare Figure 2.) 
Clearly this schema $S$ satisfies the conditions of Theorem D. 
In fact, we realize
this schema $S$ as a Expanding Hubbard Tree as follows. 
\smallskip
We join every $v_i$ in the cycle to a `fixed vertex' $p$ following that induced order from $S$ 
and making angles of $1/m$ between consecutive edges. 
(Recall that $m \ge 2$ here.) 
Now Lemma 3.2 give us a hint of how to include the vertices $\omega_0$ and $\omega_1$. 
In fact, according to Lemma 3.2, 
in the realization of this schema we must be able to find 
between $v_{j-1}$ and $\omega_1$ 
a critical point ($\omega_0$?) 
as well as between $v_{0}$ and $\omega_0$ ($\omega_1$?). 
The easiest way to satisfy these two conditions at the same time, 
is to take the ordered segment $[v_0,p]$ 
and interpolate $\omega_1$ and $\omega_0$ in that order, 
making angles of $1/d(\omega_i)$ between branches at $\omega_i$ 
(compare Figure 2). 
Defining $p$ to be fixed and of degree $1$, 
we have that with that dynamics and degree induced from $S$, this is an expanding Hubbard Tree and clearly realizes $S$. 
\medskip
\centerline{\psfig{figure=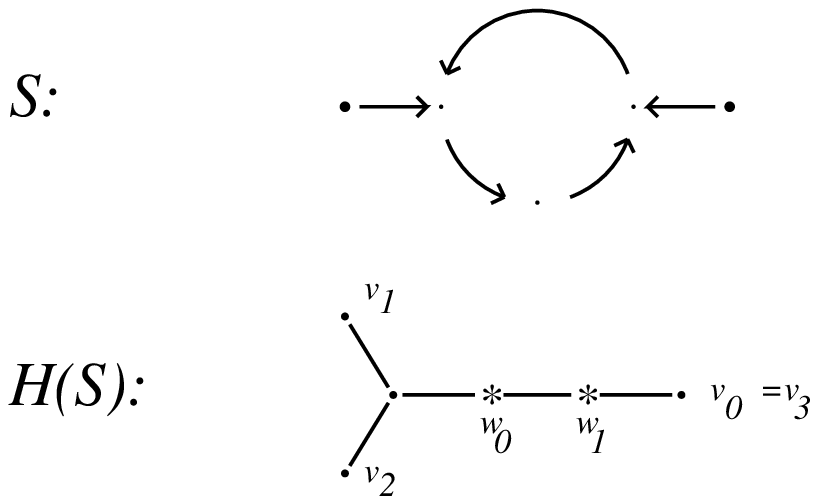,height=2.4in}} 
\noindent{\bf\sl Figure 2. Two different critical points map into a cycle of period $3$. 
Here\break 
$F(w_0)=v_1$ and $F(w_1)=v_0$. 
Note that this looks very muck like the folding construction in Lemma 8.4, Case 1.} 
\bigskip
{\bf Remark.} 
This last construction exemplifies accurately that principle we will follow 
for the realization of general schemata. 
{\it 
Whenever two points in the same fibre $W(u)$ map to the same vertex, 
there should be a critical point in $W(u)$ willing to help.} 
This principle is what motivates the construction of pseudo-chains 
(compare Section 5 and Proposition 8.8). 
\bigskip
Now we are ready to define our strategy for the proof of Theorem B. 
As we are dealing with a purely combinatorial object, 
it is natural to find first `partial solutions'. 
For example we may concentrate first on 
the cyclic part of the ambient schema ${\cal M}$; 
and in fact, Sections 4 through 8 are dedicated exclusively to this case. 
Next, even if we assume the ambient schema be cyclic, 
it is clear from the discussion and examples above that there are schemata easier to realize than others.  
In this way we will try to identify the easiest cases of ``subschemata" 
which are always present 
(compare Corollary 8.9).   
Once the easiest cases are identified and realized, 
we should learn how to `glue' partial solutions and to `append extra vertices'. 
This line of action is what justifies in part the following definition. 
\bigskip
{\bf 3.4 Admissible Subschemata.} 
Let $S$ be admissible on ${\cal M}$. 
(Here ${\cal M}$ need not be cyclic.) 
Sometimes is more appropriate to retain from ${\cal M}$ only the dynamics $F_{\cal M}$ 
(compare Definition 1.4). 
If we do so, we can define in a natural way the concept of a 
{\it consistent subschema} as follows. 
Let $S'$ be a subschema of $S$; i.e, 
the set $|S'|$ is invariant under $F$. 
Then we have an induced (ambient) subschema ${\cal M}' \subset {\cal M}$ with 
$|{\cal M}'|=\phi(|S'|)$, that dynamics induced by ${\cal M}$, and degree
$w'_{\cal M}(u)=\sum_{\{v \in |S'|:\phi(v)=u\}}w_S(v)$. 
\bigskip
{\bf Definition.} 
We say that $S'$ {\it is an admissible subschema of $S$,} 
if 
\smallskip
a) Every cycle of ${\cal M}'$ contains a critical vertex; and 
\smallskip
b) The restriction of $\phi$ to $S'$ defines an admissible projection from $S'$ to ${\cal M}'$. 
\bigskip
\eject
The actual proof of Theorem B is given in Section 9. 
However, the hardest part of the work is to realize 
that part of the schema which projects onto the cyclic part of ${\cal M}$. 
Because of this, we should concentrate most of our effort to semi-reduced schemata.  
Our proof in Section 8 will show by contradiction that every semi-reduced schema 
(admissible on a cyclic ambient schema) 
can be realized as a Hubbard Forest belonging to a special class of graphs 
to be defined in Section 6. 
(This class of Hubbard Forest has the property that 
any Hubbard Forest in this class is easy to extend, 
say by appending new vertices or by `grafting' different partial solutions.)  
However, the proof in Section 8 provides an effective recipe to actually 
extend a wrongly assumed `maximal solution' within said class. 
Therefore our proof can be taken as constructive as we provide an effective algorithm for the construction of the realization. 
The delicate point is how to start the inductive process, because the methods in Section 8 
do not always work for degree 2. 
This is related to the fact that $N(2,2)=1$ and will be treated independently in Section 7. 
\bigskip
\bigskip
\centerline{\bf 4. Critical Cycles.} 
\bigskip
\bigskip 
Suppose the ambient schema ${\cal M}$ is a cycle of period $r$ 
and $S$ consists only of critical cycles ${\cal C}_1,\dots,{\cal C}_m$. 
If this is the case, 
the period of every cycle ${\cal C}_i$ should be a multiple of $r$. 
Furthermore, given any cycle ${\cal C}_i$, 
the number of elements of this cycle 
which belong to a given fibre $W(u)$ over $|{\cal M}|$ is independent of 
the fibre and will be denoted in this section by $n(i)$. 
(This number is the return period of ${\cal C}_i$ as defined in $\S 1.7$.)
To associate a \hf 
${\bf H^*}(S)$ to $S$ we will 
construct a Forest for each cycle and then glue them together. 
For this, we will write the cycle ${\cal M}$ as 
$u_0 \mapsto u_1 \mapsto \dots \mapsto u_r=u_0$, and 
each cycle ${\cal C}_i$ as 
$v_{i0} \mapsto \dots \mapsto v_{ik}=v_{i0}$.  
We must distinguish between the cases $n(i)=1$ and $n(i)>1$. 
\medskip
For $n(i)=1$, let ${\bf H^*}(i)$ be the Forest for which 
each of its components $H_{u_j}$ 
is an edge joining points $p_{u_j}$ to $v_{ij} \in W({u_j})$. 
The degree at all $p_u$ by definition will be one, 
and they will map to each other following that order of ${\cal M}$. 
For the $v_{ij}$ the degree and dynamics is that induced by the schema $S$. 
(Here $p_u$ is introduced for two reasons; 
first, we might need later a good `gluing point' (compare Proposition 6.6), 
and second, we do not want to worry about the special cases in which the solution is a $0$-dimensional complex. 
This points $p_u$ can be dynamically interpreted as the landing point of the 
rays of argument $0$ in the realization of the forest.)
\medskip
Otherwise, if $n(i)>1$, 
we join every point in ${\cal C}_i \cap W(u)$ to a new vertex $q^i_{u}$ 
at which incident edges should form angles of $1/n(i)$ following that order induced by ${\cal C}_i$. 
We call each of these configurations a {\it star ${\cal S}_i(u)$}. 
In addition we would like to have in each component an extra Julia vertex $p_u$ as in the previous case. 
(Note that as the forest can be realized, we can always append to the forrest
the landing point of the external rays of argument $0$; 
however, we prefer to give here an explicit construction.)  
For this, we inductively associate a vertex $\omega_u \in {\cal C}_i \cap W(u)$ 
and an integer $\delta(u)>1$ to every $u \in |{\cal M}|$ in the ambient schema. 
(This will be backwards induction!) 
Let $v \in {\cal C}_i$ be any critical point in the cycle. 
Then for some $u \in |{\cal M}|$, we must have $v \in W(u)$. 
Define $\omega_u=v$ and $\delta(u)=d(\omega_u)$. 
Now suppose $F_{\cal M}(u)$ has already associated a vertex $\omega_{F_{\cal M}(u)}$ 
and an integer $\delta(F_{\cal M}(u))$. 
If $u$ has no such elements associated, 
we define $\omega_u$ as the predecessor of $\omega_{F_{\cal M}(u)}$ in ${\cal C}_i$ 
and $\delta(u)=d(\omega_u)\times\delta(F_{\cal M}(u))$. 
When all $u \in |{\cal M}|$ have these two elements associated, 
we join a vertex $p_u$ to $\omega_u$ forming an angle of $1/\delta(u)$ with the star ${\cal S}_i(u)$ (compare figure 3).
Again, the degree at all $p_u$ and $q^i_u$ should be one, 
and they map to each other in the obvious way. 
For the $v_{ij}$ the degree and dynamics is that induced by the schema. 
By construction this is clearly an expanding \hf  
${\bf H^*}(i)$. 
\smallskip 
\centerline{\psfig{figure=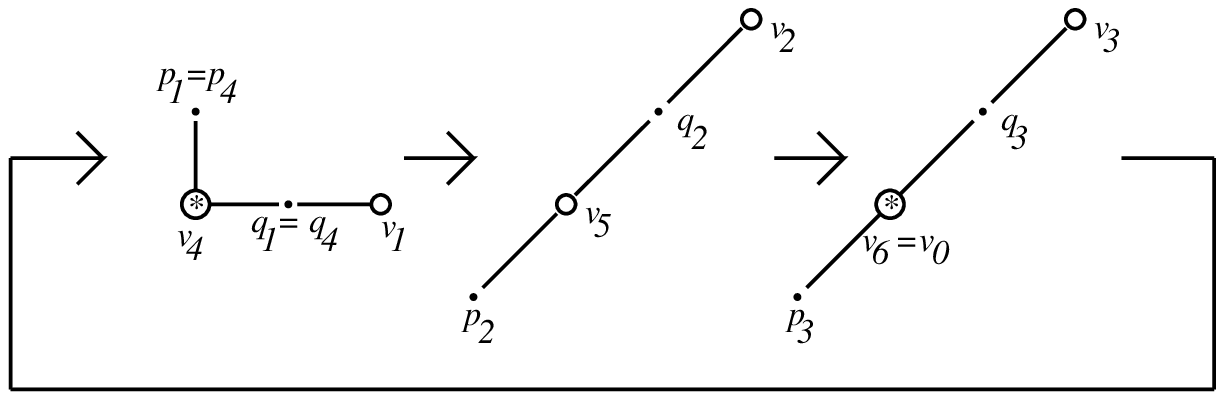,height=1.5in}} 
\noindent
{\bf\sl Figure 3. A critical cycle of period $6$ on an ambient schema of period $3$. 
Here only $v_0$ and $v_4$ are critical; at those two points the local degree is $2$. 
Note that the dynamics at $v_4$ `unbends' the tree, 
while at $v_0$ it `folds'.} 
\medskip
Finally, as there are $m \ge 1$ cycles; 
we identify the $p_u$ from different cycles, 
forming angles of $1/m$ between branches. 
This \hf  
${\bf H^*}(S)$ realizes $S$. 
\bigskip
\bigskip
\centerline{\bf 5. Pseudo-chains.} 
\bigskip
\bigskip
{\bf Definition.} 
By a {\it link} ${\cal L}(v,n)$ of lenght $n>0$ will be meant 
the ordered set\break 
$v, F_S(v), \dots, F^{\circ n}_S(v)$ formed by an element $v \in |S|$ 
and its first $n$ iterates. 
We require these $n+1$ elements to be diferrent. 
Note that by definition a link should have more than one element. 
In this case $v$ is called {\it the generator of ${\cal L}$} and 
$F^{\circ n}_S(v)$ {\it the last element}. 
The set of elements in a given link ${\cal L}$ will be denoted by $|{\cal L}|$.
\medskip
{\bf Definition.} 
Let $V \subset |S|$ be an invariant set of vertices. 
(At this stage the reader may think of $V$ as the empty set; 
compare the example below and the remark following Example 3.3.) 
An ordered collection 
${\cal L}_1={\cal L}(\omega_1,n_1), \dots, {\cal L}_k={\cal L}(\omega_k,n_k)$ 
of disjoint links and disjoint from $V$ is called a {\it V-pseudochain} if 
\smallskip
i) The set $V \cup |{\cal L}_1| \cup \dots \cup |{\cal L}_k|$ 
is invariant. 
\smallskip
ii) All generators $\omega_1, \dots \omega_r$ are critical. 
\smallskip
iii) For any $j < k$, the last element of ${\cal L}_j$ and the generator of 
${\cal L}_{j+1}$ belong to the same fibre $W(u)$ over ${\cal M}$ but 
have different images under $F_S$. 
\medskip
If all these conditions are met, we will say that 
the ordered set  $\omega_1, \dots, \omega_r$ generates the pseudo-chain. 
\smallskip
In other words, this definition simply states that `the chain 
of iteration can only be broken' 
by passing from a given vertex to a critical vertex in the same fibre. 
We require then the image of these two points be different. 
\medskip
{\bf Closed Pseudochains.} 
A $V$-pseudochain ${\cal L}_1,\dots, {\cal L}_k$ is
{\it closed} if the last element of ${\cal L}_k$ and 
the generator of ${\cal L}_1$ 
belong to same fibre over ${\cal M}$ but have different images under iteration. 
Note that in this case the order of the generating sequence of 
critical vertices $\omega_1,\dots,\omega_k$ can be cyclicly permuted. 
It follows easily from the definition that closed pseudo-chains will `project' onto a cycle within ${\cal M}$. 
\bigskip
{\bf Realizing Closed ($\emptyset$-) Pseudo Chains.} 
Suppose the elements of the admissible schema $S$ 
can be ordered so that they form a closed ($\emptyset$)-pseudochain with generators $\{\omega_1,\dots,\omega_r\}$. 
In this case we have an induced cyclic order for the vertices in $|S|$. 
For any $u \in |{\cal M}|$ we consider the set 
$W_0(u)=W(u)-\{\omega_1,\dots,\omega_k\}$. 
\medskip
{\bf Lemma.} 
{\it The number of elements in each set $W_0(u)$ is independent of $u$. 
If the schema ${\cal S}$ is not a disjoint union of critical cycles, 
then this number is strictly greater\break 
than 1.} 
\medskip
{\bf Proof.} 
This follows easily from the definitions.  
In fact, 
by simply omiting the generators $\omega_1,\dots \omega_k$ in the induced 
cyclic order descrived above, 
this new induced order ``semiconjugates" to the cyclic order in ${\cal M}$ 
which is cyclic as remarked above. 
This proves that that the number of elements in $W_0(u)$ is 
independent of $u \in |{\cal M}|$. 
\smallskip
Now, suppose $S$ is not a disjoint union of cycles. 
In this case there exists $v_1,v_2$ in the same fibre $W(u)$ such that 
$F(v_1)=F(v_2)$. 
It follows from Lemma 3.1 and the consistency of $S$ that there exists 
$v_3 \in W(u)-\{v_1,v_2\}$. 
Now, if among these three elements we can find two which are not 
generators of the pseudo-chain 
we have that $W_0(u)$ has at least two elements and we are done. 
Otherwise, if two of these are generators of the pseudochain, 
we have by definition that their respective images are not generators. 
As these images are again by definition distinct and belong to the same fibre 
$W(F_{\cal M}(u))$; 
we conclude  that $W_0(F_{\cal M}(u))$ has at least two elements. 
\endofproof
\medskip 
{\bf Realizing ($\emptyset$-) Pseudo Chains (Continue).} 
Whenever the cardinality of the set $W_0(u)$ is one, the above lemma shows 
that $S$ is the disjoint union of critical cycles. 
In this particular case $S$ can be realized as show in Section 4. 
\smallskip
Otherwise, let $m$  the cardinality of a set $W_0(u)$ be greater than 1. 
To realize $S$ as a forest we make use of this and the fact that the sets 
$W_0(u)$ have a natural cyclic order induced from the definition of 
pseudo-chain. 
We join every point in $W_0(u)$ to a new vertex $q_u$ 
following that induced cyclic order in the definition of pseudochain.  
Consecutive edges should form angles of $1/m$. 
We include the other vertices $\omega_1, \dots \omega_r$ as follows. 
To simplify notation we use the convention $\omega_0=\omega_r$. 
Now let ${\cal L}(\omega_i)$ be a link in the pseudo-chain. 
By construction $v_i$, the last point of this link is an end in the graph 
so far constructed. 
Let $u = \phi(v_i)$ (i.e, $v_i \in W_0(u)$), 
by definition of closed pseudo chain it follows that $\omega_{i+1} \in W(u)$. 
In this case we interpolate $\omega_{i+1}$ between $q_u$ and $v_i$ 
making an angle of $1/d(\omega_{i+1})$ between branches. 
The degree at all $q_u$ by definition will be one, 
and they map to each other following that order of ${\cal M}$. 
For the other vertices the degree and dynamics is that induced by the schema $S$. 
By construction this is clearly an expanding \Hf. 
\bigskip
\eject
{\bf Example.} 
Consider the schema with diagram 
\smallskip
\centerline{\psfig{figure=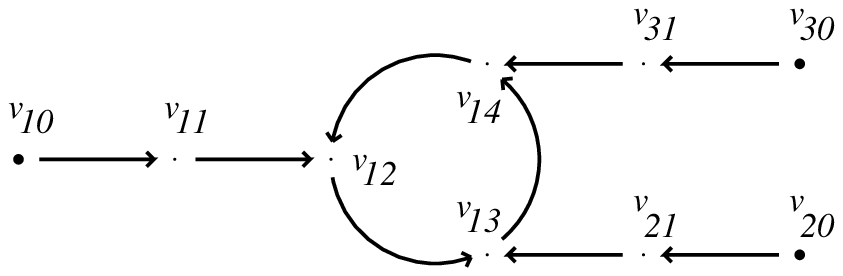,height=1.5in}} 
\smallskip
\noindent
and projection map given by $\phi(v_{ij})=i+j$ {\it (mod 3)}. 
A ($\emptyset$)-closed pseudo-chain is given by 
$$
{\bf v_{10}},v_{11},v_{12},v_{13},v_{14},{\bf v_{20}},v_{21},{\bf v_{30}},v_{31}.
$$
(Here bold face elements indicate the generating sequence of critical vertices.) 
Thus, according to the algorithm above, this schema can be realized as shown in Figure 4. 
\smallskip
\centerline{\psfig{figure=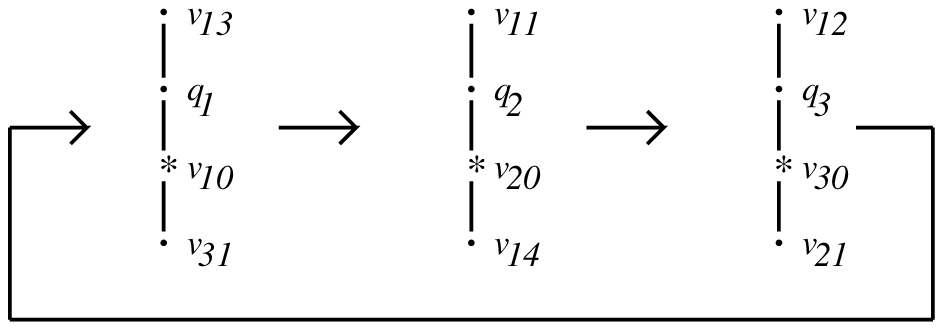,height=1.9in}} 
\noindent{\bf\sl Figure 4. Note that the relative position of $v_{31},v_{14}$ and $v_{21}$ respect to the critical points guarantees the angle condition is satisfied.} 
\bigskip
{\bf Remark.} 
Note that the algorithm described above 
guarantees that at any vertex $q_u$ in the construction, 
the `return dynamics to this component' is such that the angles are rotated $1/m$. 
Note that by a simple modification of the argument we can always achieve that 
close to this point $q_u$ the dynamics of the `return map' is rotation 
by $p/m$ whenever $p$ and $m$ are relatively prime. 
\smallskip
{\bf Remark.} 
Note that even if the schema $S$ consists of a disjoint union of critical 
cycles 
which can be aranged as a closed pseudo chain, 
the algorithm above can still be applied. 
In fact, in order to apply the algorithm we only need that every fibre 
$W_u$ contains at least two elements which are not generators. 
\smallskip
In fact, given a `limb' in the Mandelbrot set corresponding to a periodic argument of period $n$, 
the unique hyperbolic component of period $n+1$ in this limb can be reconstructed by combining the last two remarks. 
\bigskip
\bigskip
\eject

\centerline{\bf 6. Constructive Methods.} 
\bigskip
\bigskip
{\sl 
In this section we assume again the ambient schema ${\cal M}$ is cyclic. 
We define the special class of \Hf s 
which will serve as the model for the realization of semi-reduced schemata. 
All examples considered so far belong to this class. 
Within this class we show how partial solutions can be `glued' and 
how the realization of a subschema can be enlarged to include 
additional vertices.} 
\bigskip 
{\bf Digression.} 
Let $f \in {\cal P}^{\cal M}$ be post-critically finite. 
Let $P$ be the first return map of $f$ to one of the copies 
$\{u\}\times {\bf C}$ (with $u \in {\cal M}$) of the complex numbers. 
The following digression associates dynamical information of $f$ 
with dynamical properties of the return map $P$. 
\smallskip
Among obvious properties we start by 
noticing that $P$ is a polynomial of degree $n({\cal M})$ (the inner degree of ${\cal M}$). 
As periodic orbits of period $k$ for $P$ are in canonical correspondence with periodic orbits of return period $k$ for $f$, 
it follows that $f$ 
has exactly $N(n({\cal M}),k)$ cycles of return period $k$ 
(compare Theorem 1.6); 
and in particular any such polynomial map $f$ has $n({\cal M})=N(n({\cal M}),1)$ cycles of return period $1$. 
Next, as $J(P)=J(f) \cap \{u\}\times {\bf C}$, 
the study of the topological structure of the Julia set of $f$ is reduced 
to the study of the dynamical properties of $P$. 
\medskip
Let $z \in J(P)$, 
then {\it the incidence number $inc(z)$ of $z$} 
is the number of components of $J(P)-\{z\}$. 
Because $P$ is post-critically finite, 
this finite number equals the number of external rays landing at $z$ 
(compare [DH] or [P1]). 
If $z$ belongs to the Fatou set we write $inc(z)=0$. 
The relation between this incidence number $inc(z)$ and 
the number of edges incident at $z$ in a Hubbard Tree of $P$
is given by the following proposition. 
Its proof is not difficult and may be found in [P1]. 
\bigskip
{\bf 6.1 Proposition.} 
{\it 
Let $P$ be post-critically finite. 
Suppose $z \in J(P)$ is periodic and $inc(z) \ge 2$. 
Let $H$ be any Hubbard Tree associated with $P$. 
Then $z \in H$ and 
$inc(z)$ equals the number of connected components of $H-\{z\}$. 
In particular, 
if $inc(z)>2$ then $z$ is vertex in $H$ and $inc(z)$ is the number of edges of the tree incident at $z$.} 
\endofproof
\bigskip
It is also known that at least one of the fixed points should have rotation number zero. 
In fact, those fixed points of zero rotation number 
are the landing points of the `fixed rays' $R_{k/(n({\cal M})-1)}$ 
($k=0,\dots n({\cal M})-2$). 
Note in the proposition, 
that if $inc(z)=1$ for some fixed point, then necessarily $z$ has zero rotation number. 
\smallskip
Let $Z_0$ be the set of those fixed points of zero rotation number. 
Then it follows from our discussion that 
$$
n({\cal M})-1=\sum_{z \in Z_0}inc(z). 
$$
\noindent 
Those fixed points (or cycles if we are working in a \Hf ) 
of zero rotation number can be found in a Hubbard Tree as follows. 
Suppose the vertex $v$ is fixed and of Julia type. 
Let $\ell$ be any edge incident at $v$. 
Then $v$ has rotation number zero if and only if $\ell$ maps to itself. 
If follows in particular that every fixed point $v$ 
with incident number $inc(v)=1$ necessarily has zero rotation number. 
\medskip
The main point behind the next definition is that it is convenient to have control 
over the behavior of the inverse images of at least one cycle of return period 1 and zero rotation number. 
\medskip
{\bf 6.2 Tame \Hf s.} 
Let ${\bf H^*}$ be a \hf with underlying cyclic schema ${\cal M}$. 
By a {\it tame cycle in ${\bf H^*}$} 
will be meant a {\bf non critical} cycle\break 
${\cal C}:p_0 \mapsto p_1 \mapsto \dots \mapsto p_r=p_0$
of return period 1 and zero rotation number 
(in particular all $p_i$ are of Julia  type) 
and satisfying the following condition imposed to ``all the other inverses of ${\cal C}$" 
(i.e, except those in ${\cal C}$ itself):  
\smallskip
{\bf (T)}
{\it Suppose $q \in F^{-1}({\cal C})-{\cal C}$ is a vertex of ${\bf H^*}$, 
then all the angles between edges incident at $q$ are multiples of $1/d(q)$.}
\smallskip
In other words, the `germs' of edges at any preimage of that cycle 
(except for those in ${\cal C}$ itself), 
all map to the same `germ' of edge at a vertex in the cycle. 
In particular, 
if the incidence number of a return 1 periodic point is $1$, 
then that periodic point generates a tame cycle because then 
condition ${\bf (T)}$ is trivially satisfied. 
Note that 
{\it ``tameness" is a property of a cycle in the combinatorial object and 
in general makes no sense to extend this concept to cycles for maps in ${\cal P}^{\cal M}$.} (Compare Figure 5.) 
\medskip
{\bf Definition.} 
A \hf with underlying cyclic schema ${\cal M}$ is {\it tame} if it contains a tame cycle. 
More generally, a  \hf with underlying cyclic schema ${\cal M}$ 
is {\it tame} if it can be extended to a tame Hubbard Forest in the sense described above. 
Again, ``tameness" is a property of the combinatorial object and {\bf not} of the map it realizes. 
In fact, given a post-critically finite map $f \in {\cal P}^{\cal M}$, 
it may happen that two different set of invariant vertices 
define respectively a tame and a non tame \hf 
which realize this given $f$. 
\smallskip 
\centerline{\psfig{figure=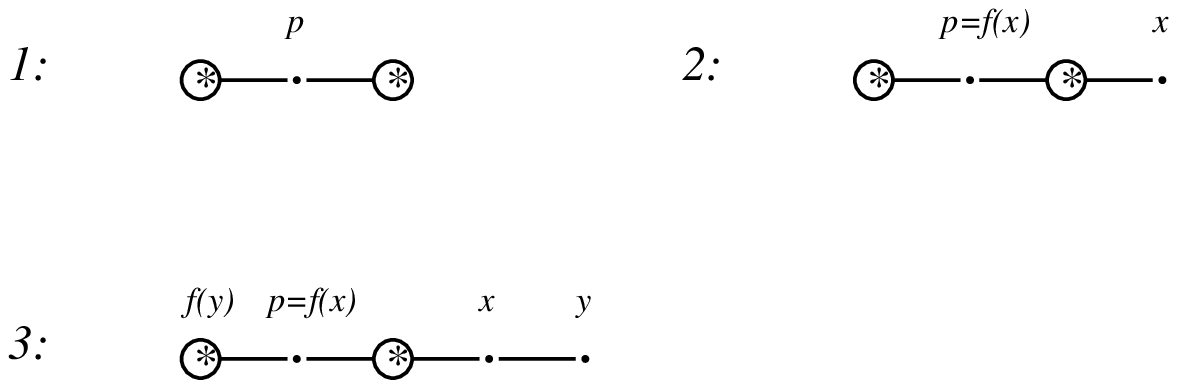,height=2.25in}} 
\smallskip
{\sl Figure 5. These three Hubbard Trees realize the dynamics of $P(z)=z^3+{3 \over 2}z$. 
Here both critical points are fixed. 
The other fixed point is between these two critical points. 
Only the first two Hubbard Trees are tame: 
In fact, in the first tree condition ({\bf T}) is vacously satisfied. 
In the second and third, 
condition ({\bf T}) should be verified only at x. 
However, at x the local degree is 1, so only the second tree is tame.}  
\eject 
In practice it is not necessary to extend a \hf 
in other to verify whether is tame. 
The following proposition gives a sufficient criterion for tameness. 
\bigskip
{\bf 6.3 Proposition.} 
{\it Let ${\bf H^*}$ be a cyclic \hf with underlying ambient schema ${\cal M}$. 
Fix any component $H_u$ of ${\bf H^*}$. 
Denote by $Z_0$ the set of those return 1 periodic vertices in ${\bf H^*}$ 
which have rotation number zero and belong to $H_u$. 
If 
$$
\sum_{v \in Z_0} inc(v) < n({\cal M})-1
$$
then ${\bf H^*}$ is tame.} 
\smallskip
{\bf Proof.} 
The proof will follow immediately from the definition and Proposition 6.1. 
Let $f \in {\cal P}^{\cal M}$ be the realization of ${\bf H^*}$. 
The hypothesis implies there is a return 1 periodic point $p$ of zero rotation number 
in that copy of the complex numbers associated with $H_u$ 
which is not a vertex in $H_u$. 
If follows from Proposition 6.1 that either $inc(p)=1$ or $inc(p)=2$. 
If $inc(p)=1$ then there is an extension of ${\bf H^*}$ which includes the orbit of $p$. 
It follows again by Proposition 6.1 that only one edge meet at $p$ in this extension. 
Therefore, this is a tame cycle. 
Otherwise, if $inc(p)=2$ then $p \in H_u$ but not as a vertex. 
It follows that an extension of ${\bf H^*}$ including the orbit of $p$ 
can be constructed without including any vertex other than those in this orbit. 
It follows that condition ${\bf (T)}$ 
should be verified over an empty set in this extension. 
\endofproof
\bigskip
Henceforth, without loss of generality we assume that whenever a \hf 
is tame, 
a tame cycle can be found among its set of vertices. 
(There is no need to reconstruct that map $f \in {\cal P}^{\cal M}$ 
which realizes ${\bf H^*}$ to obtain such extension. 
In fact, an algorithm to extend a Hubbard Tree to include the landing point  
of any ray is given in [P3, Proposition III.4.5].) 
\medskip
{\bf Grafting.} 
The importance of this class of \hf 
is given by the fact that two of these objects projecting to the same cyclic 
ambient space can be ``grafted" along their tame cycles. 
In fact, let ${\bf H_1}$ and ${\bf H_2}$ be expanding Hubbard Forests 
with tame cycles ${\cal C}_1$ and ${\cal C}_2$ respectively. 
Let ${\cal C}_i:q_{i0} \mapsto q_{i1} \mapsto \dots \mapsto q_{ir}=q_{i0}$ denote those cycles. 
We identify $q_{0j}$ with $q_{1j}$ and give the structure of a Hubbard Forest to this new object ${\bf H^*}$ as follows. 
Let $m_1$ and $m_2$ be the number of incident edges at any vertex of the cycles 
${\cal C}_1$ and ${\cal C}_2$ respectively. 
Note that these numbers are independent of that chosen vertices. 
To give ${\bf H^*}$ the structure of an expanding \hf 
it is enough to define the angles between edges at any vertex of the cycle ${\cal C}_1 \sim {\cal C}_2$ to be a non trivial multiple of $1/(m_1+m_2)$ and 
`pull back' the definition along this cycle in a compatible way.  
In fact, the expanding condition between adjacent periodic Julia vertices
and 
the angle condition at every vertex not mapping into ${\cal C}_1 \sim {\cal C}_2$  are clearly satisfy. 
However, at a vertex mapping into ${\cal C}_1 \sim {\cal C}_2$ 
condition ${\bf (T)}$ is  satisfied, 
which proves not only that ${\bf H^*}$ is an expanding Hubbard Forest, 
but is also tame.  
In fact, this `same' cycle ${\cal C}_1 \sim {\cal C}_2$ is still tame. 
\bigskip
{\bf Remark.} 
Note that as a constructive method, 
this procedure is well defined even without the assuption of tameness. 
In fact, in the general case we will need only to inductively 
redefine the angles between edges at iterated inverses of the cycle.  
The angle conditions are trivially satisfied at all other vertices 
because the local dynamics is copied from the original structures. 
Furthermore, the expanding condition between adjacent periodic 
vertices is also trivially satisfied. 
In fact, any two adjacent periodic vertices in the new forest are also 
adjacent in the original ones; 
as the images in one forest are symbolically the same as in the old forests, 
the expansiveness at some iterate follows. 
\bigskip
{\bf 6.4 Definitions.} 
Let $S$ be an admissible schema on the cyclic ambient ${\cal M}$. 
This schema $S$ is {\it tame}, 
if it can be realized by a tame \Hf. 
\smallskip
A cycle ${\cal C} \subset S$ is {\it superfluous}, 
if the connected component of $S$ which contains ${\cal C}$ has no critical vertices. 
The technical goal is to prove that 
every semi-reduced schema is tame. 
(In particular there should be no superfluous cycles in $S$.) 
So far, we have only met with tame schemata: 
\bigskip
{\bf 6.5 Lemma.} 
{\it The following schemata are tame: 
\smallskip
i) A disjoint union of critical cycles ($\S$4). 
\smallskip
ii) A closed ($\emptyset$-) pseudo-chain ($\S$5). 
\smallskip
iii) Example $\S$3.3.}
\medskip
{\bf Proof.} 
In fact, in each case we provided a tame realization. 
\endofproof
\bigskip
\bigskip
{\bf 6.6 Proposition.} 
{\it Let $S_1$ and $S_2$ be disjoint tame schemata projecting to the same cyclic set ${\cal M}$. 
Suppose neither $S_1$ nor $S_2$ have superfluous cycles. 
Then the disjoint union $S_1 + S_2$ is tame.} 
\medskip
The proof is based in the following technical lemma (compare also Lemma 8.1). 
\bigskip
{\bf 6.7 Lemma.} 
{\it 
With the hypothesis above, 
suppose ${\bf H^*}(S_1)$ is a tame realization of $S_1$. 
Then there is an extension ${\bf H^{'*}}(S_1)$ 
for which a return 1 cycle is present but is not identified with any $v \in |S_1|$ (i.e, is superfluous). 
In particular this cycle is disjoint from the orbit of the critical set. 
Furthermore, this extension can be chosen to be tame.} 
\medskip
{\bf Proof.} 
Suppose $S_1$ has $m$ components; 
because each of these components contains a critical point 
it follows easily that $n({\cal M})$ is bigger than $m$. 
This means that among the $n({\cal M})$ return 1 periodic cycles 
of the realization of ${\bf H^*}(S_1)$ as a map $f \in {\cal P}^{\cal M}$, 
at least one is disjoint from the orbit of the critical set. 
Now, if ${\bf H^*}(S_1)$ contains all these 
$n({\cal M})$ return 1 periodic cycles 
it follows that ${\bf H^*}(S_1)$ itself has the required properties. 
Otherwise, suppose there is a point $p$ which generates a period 1 cycle for $f$ 
which is not present as a vertex in ${\bf H^*}$. 
If $inc(p)=1$, we include the orbit of $p$ in the forest. 
As remarked several times above, this cycle is tame. 
If $inc(p)=2$ then the inclusion of the orbit of $p$ in the forest amounts only to 
the addition of a periodic set of vertices (compare Proposition 6.1). 
In this second case the original tame cycle will remain tame after extension. 
\endofproof
\bigskip
{\bf Proof of Proposition 6.6.} 
Let ${\bf H^*}(S_i)$ ($i=1,2$) 
be the realizations guaranteed by the hypothesis. 
Then there are tame cycles ${\cal C}_i:p_{i0} \mapsto p_{i1} \mapsto \dots \mapsto p_{ir}=p_{i0}$ which we assume to be present in  ${\bf H^*}(S_i)$. 
The grafting procedure defined above has a minor drawback for the realization of $S_1 + S_2$: 
it may happen that both tame cycles ${\cal C}_i$ are post-critical (compare Figure 6). 
If this is the case, note that the grafted Hubbard Forest will not realize 
the {\bf disjoint} union of $S_1$ and $S_2$. 
In fact, the grand orbit of the cycle will intersect in a non trivial way 
both $|S_1|$ and $|S_2|$. 
However this is impossible because those sets are by hypothesis invariant and disjoint. 
\smallskip
In order to deal with this situation we must redefine the image of the inverses of at least one of the gluing cycles. 
Here is where the tameness condition of the forest plays a central role:  
According to Lemma 6.7 we may assume there is a return 1 cycle 
$q_{10} \mapsto q_{11} \mapsto \dots \mapsto q_{1r}=q_{10}$ 
in ${\bf H^*}(S_1)$, 
which is disjoint from the orbit of the critical set 
(and therefore does not realizes any cycle in $S_1$ because superfluous cycles are not allowed by hypothesis.)  
With this in mind, we graft ${\bf H^*}(S_1)$ and  ${\bf H^*}(S_2)$ 
along the tame cycles ${\cal C}_i$ to construct a graph ${\bf H^*}$. 
The angles are redefined along `the cycle' in the obvious way described above. 
However, the dynamics at every vertex 
$v \in F_{S_2}^{-1}({\cal C}_2)-{\cal C}_2$ should be modified as follows. 
Suppose $v \in F_{S_2}^{-1}({\cal C}_2)-{\cal C}_2$ is such that 
$F_{S_2}(v)=p_{2j}$ for some $p_{2j}$ in the tame cycle ${\cal C}_2$. 
In this case (and only in this case) we redefine the dynamics at $v$ as 
$F_{\bf H^*}(v)=q_{1j}$. 
In other words, we  ``push" the image of $v$ from the tame cycle ${\cal C}_2$ to the superfluous cycle in ${\bf  H^*}(S_1)$ guaranteed by Lemma 6.7. 
All other structure is defined in ${\bf H^*}$ as it was defined for 
${\bf H^*}(S_1)$ and ${\bf H^*}(S_2)$ 
(except of course, for the angles along the gluing cycle).  
Clearly this new combinatorial object ``realizes'' $S=S_1 + S_2$. 
In fact, if a cycle in $S_2$ is realized in ${\bf H^*}(S_2)$ by ${\cal C}_2$, 
then this cycle is realized in ${\bf H^*}$ by 
$q_{10} \mapsto q_{11} \dots q_{1r}=q_{10}$. 
All other cycles in $S_1 + S_2$ are realized as they were before. 
(Compare Figure 6.) 
To complete the proof, 
it is enough to prove that we have an expanding Hubbard Forest. 
\medskip
By construction, every connected component of ${\bf H^*}$ is a tree. 
Also a vertex in ${\bf H^*}$ is of Julia type if and only if 
the corresponding vertex in either ${\bf  H^*}(S_1)$ or ${\bf  H^*}(S_2)$ is so. 
From this it follows easily that 
the expanding condition between adjacent periodic Julia vertices in ${\bf H^*}$, as well the normalization for angles at Julia vertices are satisfied. 
\smallskip
\centerline{\psfig{figure=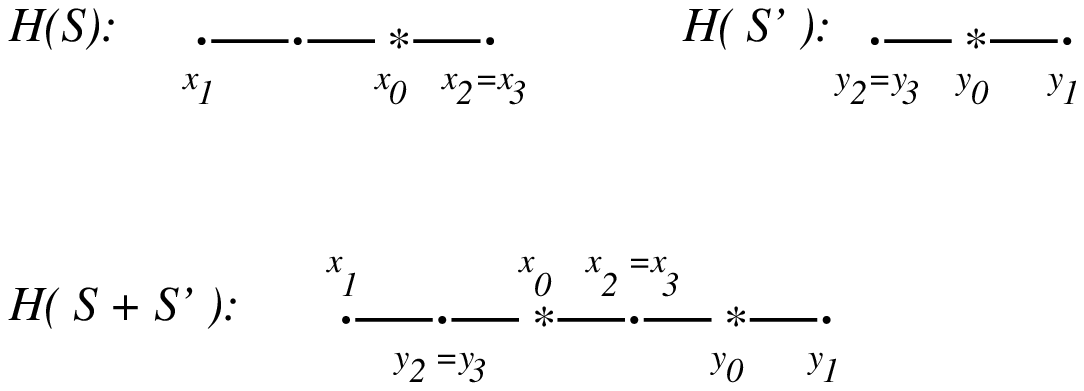,height=2.5in}} 
\smallskip
\noindent{\sl Figure 6. The image of $y_1$ is pushed from the tame cycle in ${\bf H(S')}$ to the superfluous cycle in ${\bf H(S)}$.}
\bigskip
\noindent
However, there are two delicate points which we still need to check. 
First, suppose $\ell$ is an edge with adjacent vertices $v,v'$ where 
$v \in F_{S_2}^{-1}({\cal C}_2)-{\cal C}_2$. 
We need to show that $\ell$ viewed as a germ of edge at $v'$, 
still maps to the same germ at $F_{S_2}(v')$. 
But this is clear because $F_{\bf H^*}(\ell)$ can be visualized as 
the disjoint union of $F_{S_2}(\ell)$ in ${\bf H^*}(S_2)$ and the arc $[p_{1j},q_{1j}]$ in ${\bf H^*}(S_1)$. 
The second delicate point is if the angle condition is satisfied at 
$v \in F_{S_2}^{-1}({\cal C}_2)-{\cal C}_2$. 
But this follows immediately from the fact that the cycle ${\cal C}_2$ is tame in  
${\bf H^*}(S_2)$. 
Finally the cycle ${\cal C}_1$ in ${\bf H^*}$ is clearly tame 
bacause it has exactly the same inverses in ${\bf H^*}$ 
as it had in ${\bf H}(S_1)$. 
\endofproof 
\bigskip
{\bf 6.8 Saturated Points.} 
Let $S$ be admissible on ${\cal M}$. 
Pick a vertex $u \in |{\cal M}|$.  
We say that a vertex $v \in W_{F_{\cal M}(u)}$ of $S$ is {\it $u$-saturated} if 
$$
\sum_{\{v' \in W(u):f(v')=v\}} d(v') = d_{\cal M}(\phi(v)).
$$
In other words the number of inverses of $v$ in a given fibre $W(u)$ is maximal counting multiplicity. 
In particular note that if ${\cal M}$ is cyclic, every end of $S$ is not saturated because it has no inverses. 
Note that in the definition there is no need to assume ${\cal M}$ cyclic. 
However, if ${\cal M}$ is not cyclic, a distinction must be made respect which fibre the saturation is referred to. 
\bigskip
The next proposition simply states that under certain conditions, 
the realization of a given schema can be enlarged by appending one preperiodic
vertex a the time in the forest 
and then ``criticalize" this vertex if necessary. 
Note that this can be done if for a given vertex not all its possible 
inverses are compromized in the realization of the starting schema. 
This fact is what is encapsuled by the definition of saturation.  
\bigskip
{\bf 6.9 Proposition.} 
{\it Let $S$ be admissible on the cyclic schema ${\cal M}$. 
Suppose that $|S|$ can be described as the disjoint union of vertices $|S'|$ 
of an admissible subschema $S'$ which can be realized as a Hubbard Forest ${\bf H^*}(S')$, 
and a portion of the orbit ${\cal O}(\omega)$ of a critical point 
$\omega \not\in |S'|$. 
Suppose further that $|S'| \cap {\cal O}(\omega)$ is not empty. 
Denote by $v$ that last point in the orbit of $\omega$ which does not belong to $|S'|$. 
If $F(v)$ is not $\phi(v)$-saturated respect to 
the induced ambient schema ${\cal M}'$ of $S'$, then $S$ is tame.}   
\medskip
{\bf Proof.} 
As $F(v) \in V_{\phi(F(v))}$ is not $\phi(v)$-saturated, 
we can find an extension of ${\bf H^*}(S')$ for which there is an inverse $v \not\in |S'|$ of $F(v)$. 
In particular this point $v$ is not critical in this extension, 
so if $v \in |S'|$ has degree $d(v)=1$ then this extension 
realizes $\{v\}+S'$. 
Otherwise, if $d(v)>1$ we modify this Hubbard Forest to realize $\{v\}+S'$ as follows. 

If $v \in |S|$ has degree $d(v)$, 
we define the local degree of $v$ in the tree to be also $d(v)$ and 
modify the angle function between edges at $v$ 
so that the angle condition is satisfied. 
(We have ``criticalized" the vertex $v$.)
That a modification compatible with the dynamics is possible 
follows from the fact that $v$ is preperiodic. 
Next, we inductively proceed to modify the angle function at every iterated inverse of $v$ which is a vertex in the original forest. 
\smallskip
Now, the same hypothesis of the proposition applies to the schema $S' + \{v\}$ 
(in which trivially $v$ is not saturated) 
and the remaining portion of the orbit ${\cal O}(\omega)$. 
But this means that after a finite number of steps we are done. 
(Remark: Note that this same proof applies if the schema $S$ has superfluous cycles. 
In fact, we have not imposed any restriction in the hypothesis beyond that 
that $S'$ can be realized; 
in particular nothing is said about the tameness of $S'$.)  
Furthermore, 
{\it this \hf ${\bf H^*}(S)$ has the same number of periodic orbits 
(and of the same period) as does ${\bf H^*}(S')$;}  
this because at each step we have only modified the angle function and the degree but not the dynamics.  
\smallskip
Finally, as the critical point $\omega$ does not belong to the schema $S'$, 
it follows that the inner degree $n({\cal M})$ of ${\cal M}$ is bigger than that $n({\cal M}')$ of the ambient schema ${\cal M}'$ of $S'$ 
(compare Definition 3.4). 
It follows easily from Proposition 6.3 that ${\bf H^*}(S)$ is tame. 
\endofproof
\bigskip 
\bigskip
\centerline{\bf 7 Subordinated Configurations.} 
\bigskip
\bigskip
{\sl Finally we consider an extra special case in which the tools developed so far may not apply. 
(This last special case is closely related to the uneasy situation where $N(2,2)=1$.) 
Again we assume the ambient schema ${\cal M}$ is cyclic.} 
\bigskip
{\bf Definition.} 
A {\it subordinated configuration $S$} 
is an admissible schema $S$ 
which is the disjoint union of two components $S_1$ and $S_2$. 
Each $S_i$ should be generated by a critical point $v_{i0}$. 
In other words, $S$ can be described as the disjoint union of the orbits of those two critical points. 
We require $S_1$ to be non admissible and $S_2$ to be tame. 
(Thus, $S_1$ is {\it subordinated to $S_2$}.) 
\bigskip
First we briefly study those conditions implied by the hypothesis 
and establish the notation to be used throughout the rest of this section. 
Let ${\cal C}_i$ be the unique cycle contained in $S_i$. 
\medskip
{\bf 7.1} 
As $S_2$ is tame (and in particular admissible) 
we have the following dichotomy. 
Either the generating critical point $v_{20}$ of $S_2$ is periodic or not. 
In the later case there are exactly two different points 
$v_{2i_2},v_{2j_2}$ in the same fibre which satisfy 
$F(v_{2i_2})=F(v_{2j_2})$. 
(Say $v_{2j_2} \in {\cal C}_2$ and $v_{2i_2} \not\in {\cal C}_2$.) 
It follows from Lemma 3.1, 
that that same fibre, say $W(u_2)$,  
contains a critical vertex $v \in S_2$ which satisfies $F(v) \ne F(v_{2j_2})$.
\smallskip
In the former case, that is, when $v_{20}$ is periodic, 
we set for consistency $v_{2j_2}$ as the unique vertex which maps to $v_{20}$. 
Thus, the orbit of $v_{20}$ will always be described as 
$$
v_{20}, v_{21}, v_{22}, \dots, v_{2j_2}. 
$$
\medskip
{\bf 7.2} 
As $S_1$ is not admissible, 
the generating critical point $v_{10}$ of $S_1$ is not periodic 
(compare $\S 4$) 
that is $v_{10} \not\in {\cal C}_1$.  
Again it follows that there are exactly two different points 
$v_{1i_1},v_{1j_1}$ in the same fibre, say $W(u_1)$, which satisfy 
$w=F(v_{1i_1})=F(v_{1j_1})$. 
(Say $v_{1j_1} \in {\cal C}_1$ and $v_{1i_1} \not\in {\cal C}_1$.) 
Because condition a) in the definition of admissible schema can only fail at $w$ 
(compare Definition 1.8) 
we necessarily have
$$
d(v_{1i_1})+d(v_{1j_1}) > 1+ \sum_{v' \in W(u_1) \cap |S_1|} (d(v')-1). 
$$ 
(Here the right hand side of the expression denotes the degree of $u_1 \in |{\cal M}|$ viewed as the ambient schema of $S_1$.) 
As $v_{1i_1},v_{1j_1} \in W(u_1)$ it follows easily that the fibre $W(u_1)$ 
contains no critical points in $S_1$ except perhaps $v_{1i_1}$ or $v_{1j_1}$.
\smallskip
The dichotomy here is between the critical point 
$v_{10}$ belonging to the fibre $W(u_1)$ or not. 
If $v_{10} \in W(u_1)$ then necessarilly $v_{10}=v_{1i_1}$ and therefore $F(v_{10})$ is periodic. 
\medskip
{\bf 7.3} 
Finally, as $S_1+S_2$ is admissible, 
it follows that the fibre $W(u_1)$ (for $S_1 + S_2$) contains a critical vertex $v$ 
other than $v_{1i_1}$ or $v_{1j_1}$ 
(compare Lemma 3.1). 
The discussion in $\S$7.2 shows that $v \in |S_2|$. 
\bigskip
{\bf 7.4} 
As Proposition 6.6 clearly can not be applied in this context, 
we turn our approach to Proposition 6.9. 
For example we will like to know under which conditions the schema ${\cal C}_1 +S_2$ can be realized. 
If this is the case, 
it will follow immediately from Proposition 6.9 that $S_1 + S_2$ is tame 
(compare Lemma 7.5). 
However, the conditions for the realization of ${\cal C}_1 +S_2$ are hardly a surprise 
(compare Lemma 7.6). 
\bigskip
{\bf 7.5 Lemma.} 
{\it Under the hypothesis and notation above, suppose 
the schema ${\cal C}_1 + S_2$ can be realized by a Hubbard Forest ${\bf H^*}$. 
Then the schema $S_1 + S_2$ is tame.} 
\medskip
{\bf Proof.} 
It follows from the discussion in $\S$7.3 that 
the vertex $w=F(v_{1i_1})$ ($=F(v_{1j_1}))$ is not 
$u_1$-saturated in the schema ${\cal C}_1+S_2$.
Therefore the hypothesis of Proposition 6.9 apply with $S'={\cal C}_1 + S_2$ 
and $\omega=v_{10}$. 
The result follows. 
\endofproof
\bigskip 
{\bf 7.6 Lemma.} 
{\it Under the hypothesis and notation above,  
the schema ${\cal C}_1 + S_2$ can be realized by a Hubbard Forest ${\bf H^*}$ 
if and only if 
${\cal C}_1 + S_2$ is admissible.} 
\medskip 
{\bf Proof.} 
If ${\cal C}_1 + S_2$ can be realized by a Hubbard Forest ${\bf H^*}$ 
then Theorem A shows that ${\cal C}_1 + S_2$ is admissible. 
Conversely, if ${\cal C}_1 + S_2$ is admissible we distinguish between whether ${\cal C}_1$ is a critical cycle or not. 
If ${\cal C}_1$ is a critical cycle, Lemma 6.5 says that ${\cal C}_1$ is also tame and the result then follows from Proposition 6.6. 
Otherwise if ${\cal C}_1$ contains no critical points, 
the realization of the tame schema $S_2$ as a map in ${\cal P}^{\cal M}$ 
contains a cycle of the same return period as ${\cal C}_1$ which is not post-critical (compare Lemma 8.1 and Theorem 1.6). 
The result follows by including this cycle in the Hubbard Forest and 
identifying it with ${\cal C}_1$. 
(Compare Figure 7.) 
\endofproof
\medskip
\centerline{\psfig{figure=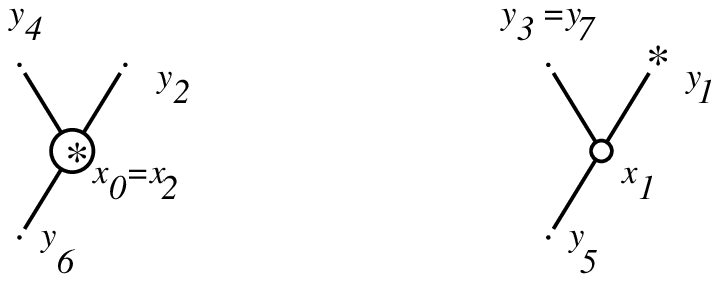,height=1.5in}} 
\smallskip
{\sl Figure 7. 
Here the orbit of $y_1$ is ``subordinated" to that of $x_0$. 
(We left to the reader the easy task of determining the angle between branches.)}
\bigskip
\bigskip
{\bf 7.7 Lemma.} 
{\it Under the hypothesis and notation above,  
suppose the schema ${\cal C}_1 + S_2$ is not admissible. 
Then $v_{20}$ is the unique critical point of $S_2$ and has degree $d(v_{20})=2$. 
Furthermore, both ${\cal C}_1$ and ${\cal C}_2$ have return period 2 and 
${\cal C}_1$ contains no critical vertices.}   
\medskip 
{\bf Proof.} 
The relation $N(n,k)<2$  
holds if and only if $n=k=2$. 
\endofproof
\bigskip
{\bf 7.8 Remark.} 
Suppose ${\cal V} \subset |S_1|$ is an invariant set of vertices properly contained in $|S_1|$ 
(in particular $v_{10} \not\in {\cal V}$). 
If ${\cal V} + S_2$ can be realized, then $S_1 + S_2$ is tame. 
In fact, in this case we can use again Proposition 6.9 as in Lemma 7.5. 
\bigskip
It follows that in order to prove that 
every subordinated configuration is tame, 
it is enough to assume that $v_{10}$ and $v_{20}$ are the only critical vertices in $S_1 + S_2$, 
and that furthermore, both ${\cal C}_1$ and ${\cal C}_2$ have return period 2 
with ${\cal C}_1$ containg no critical vertices. 
In fact, according to our previous discussion 
all other cases are considered in Lemmas 7.4-7 or can be reduced to this one using Remark 7.8. 
Therefore, it is enough to consider that cases where the cyclic ambient schema ${\cal M}$ 
contains one or two points, corresponding to the cases whether $v_{10}$ and $v_{20}$ belong to the same fibre or not. 
In what follows we freely use the notation established in $\S$7.1-4. 
\bigskip
{\bf 7.9.} 
Suppose that $v_{10}$ and $v_{20}$ belong to the same fibre $W(\phi(v_{10}))=W(\phi(v_{20}))$. 
In this case we can assume without loss of generality 
that $|{\cal M}|=\{u_1\}$.  
First suppose $v_{20}$ does not belong to a critical cycle. 
In this case 
$$
{\bf v_{20}}, v_{21}, \dots, v_{2j_2}, {\bf v_{10}}, v_{11}, \dots, v_{1j_1} 
$$ 
is a closed ($\emptyset$-) pseudo chain with $\{v_{20},v_{10}\}$ as generators, 
and therefore $S_1 + S_2$ is tame by Lemma 6.5.  
In fact $v_{2j_2}$ and $v_{10}$ 
(as well as $v_{1j_1}$ and $v_{20}$) 
by definition belong to different components of $S_1 + S_2$ and therefore 
map to different points. 
\smallskip
Otherwise, assume $v_{20}$ belongs to a critical cycle. 
It follows from the discussion in $\S$7.2 that $v_{11}=F(v_{10})$ is periodic. 
Thus in this case $S_1 + S_2$ is tame as is shown in Figure 8. 
\smallskip
\centerline{\psfig{figure=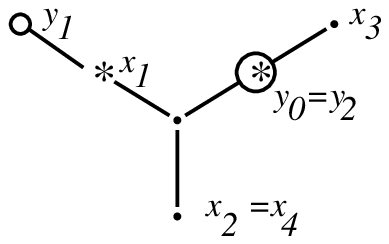,height=1.8in}} 
\smallskip
\noindent{\sl Figure 8. A tame Hubbard Tree. Here the orbit of $x_1$ is subordinated to that of $y_0$.}
\bigskip
\bigskip
{\bf 7.10.} 
Suppose now that $v_{10}$ and $v_{20}$ belong to different fibres. 
We can assume then that $|{\cal M}|$ has exactly two points. 
From our preliminary discussion, 
it follows that the elements of $S_1 + S_2$ can be ordered as a pseudo-chain 
$$
{\bf v_{10}},v_{11},\dots,v_{1j_{1}},{\bf v_{20}},v_{21},\dots,v_{2j_{2}}.
$$
In fact, from the discussion in $\S$7.3 follows that $v_{20} \in W(u_1)$; 
by definition (compare $\S$7.2) we have $v_{1j_{1}} \in W(u_1)$, 
Finally $F(v_{20}) \ne F(v_{1j_1})$ as these vertices belong to different components of $S_1+S_2$. 
\bigskip
{\bf 7.11.} 
First we study that case where $v_{20}$ is periodic. 
In this case, $v_{2j_2}$ was defined as the unique preimage of $v_{20}$. 
It follows that $v_{20}$ and $v_{2j_2}$ belong to different fibres. 
Therefore, $v_{10}$ and $v_{2j_2}$ belong to the same fibre and  
have different images as they belong\break 
to disjoint components of $S_1 + S_2$. 
Therefore that pseudo-chain defined in $\S$7.10 is closed. 
But then $S_1 + S_2$ is tame as was shown in Section 5 (compare also Lemma 6.5). 
\bigskip 
{\bf 7.12.} 
If $v_{20}$ is not periodic, 
that pseudo-chain is not closed. 
In fact, 
the discussion in $\S$7.1 shows that $v_{2j_2}$ shares the same fibre with a critical point belonging to $|S_2|$. 
As $v_{20}$ is by assumption that unique critical point, the claim follows. 
Moreover, 
that fibres $W_1=W(\phi(v_{10}))$ and $W_2=W(\phi(v_{20}))$ to which the critical points 
$v_{10}$ and $v_{20}$ respectively belong have different cardinalities. 
Following that same ideas as in Section 5 for the realization of pseudo-chains, 
we define sets 
${\hat W}_1=W_1-\{v_{10}\}$ and 
${\hat W}_2=W_2-\{v_{20},v_{2j_2}\}$. 
It follows that these two sets have the same cardinality $m \ge 2$. 
(This last claim follows from the fact that ${\cal C}_1$ has 4 elements by assumption.) 
As in Section 5, we join every point in 
${\hat W}_i$ to a vertex $p_i$ following that induced order in the pseudochain, 
and forming angles of $1/m$ between consecutive edges. 
Now $v_{2j_2-1}$ is the last element in the pseudo-chain that shares that same fibre with $v_{10}$. 
We insert $v_{10}$ between $v_{2j_2-1}$ and $p_1$ forming angles of $1/d(v_{10})$.  
Similarly $v_{20}$ is interpolated between $v_{1j_1}$ and $p_2$ forming angles of $1/d(v_{20})$. 
\smallskip
Three important remarks will show us how to include the missing vertex $v_{2j_2}$ 
(compare Figure 9). 
First note that because the critical vertex $v_{10}$ is between $v_{2j_2-1}$ and $p_1$, 
and $v_{11}$ is an `end' of the so far constructed graphs, 
$v_{2j_2}$ might be (almost) anywhere in the tree 
(except in that segment joining $p_2$ to $v_{11}$). 
Second note that $F(v_{2j_2})$ has two preimages from which we can so far find one in the graph. 
(Namely $v_{2i_2}$ which is clearly different from $v_{1j_1}$.) 
According to Lemma 3.4, we necessarily conclude that $v_{2j_2}$ 
should {\bf not} be included in the `same side of $v_{20}$ as $v_{2i_2}$'. 
As $v_{20}$ was interpolated between $v_{1j_1}$ and $p_2$, it follows that 
$v_{2j_2}$ should be included in that `same side of $v_{20}$ as $v_{1j_1}$'. 
Finally, if this is to be the case, 
the regulated set containing $v_{20},v_{1j_1},v_{2j_2}$ should be homeomorphic to 
the regulated set containing $F(v_{20}),F(v_{1j_1}),F(v_{2j_2})$ and with the same angles between corresponding edges. 
But these last three points are already included in the same component of our graph. 
Therefore the tree they generate can serve as model for that generated by $v_{20},v_{1j_1},v_{2j_2}$. 
Thus, the best way to realize $S$ is to delete the segment $[v_{20},v_{1j_1}]$, 
and after identifying $F(v)$ with $v$ replace it by a copy of the 
angled tree generated by $F(v_{20}),F(v_{1j_1}),F(v_{2j_2})$. 
The angle at $v_{20}$ should be $1/d(v_{20})$. 
This is clearly a tame \hf because there is only one return 1 cycle present; 
namely, 
that orbit of $p_1$ and $p_2$ (compare Proposition 6.3). 
\endofproof
\smallskip
\centerline{\psfig{figure=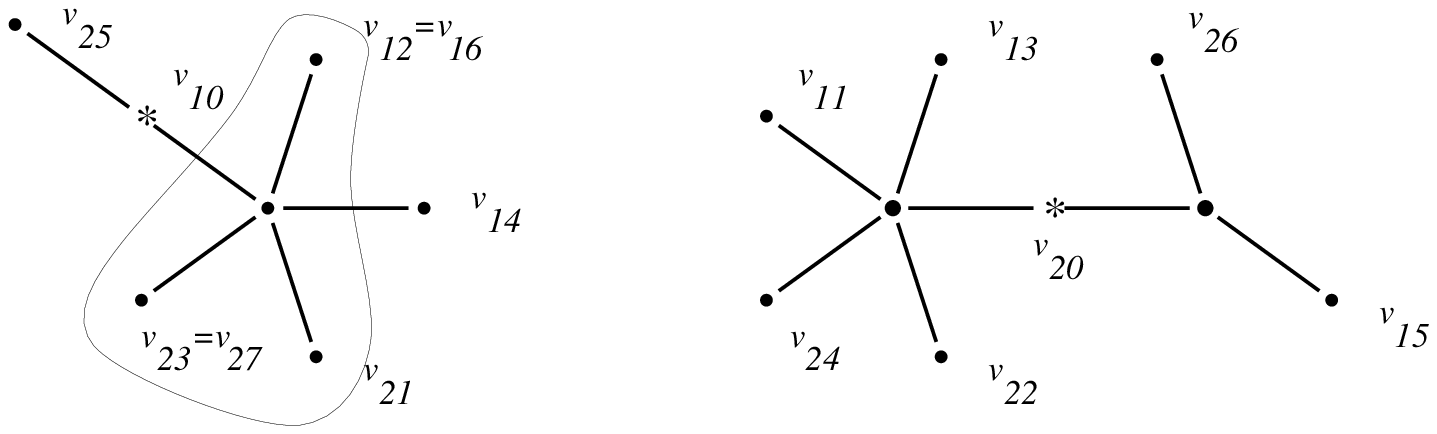,height=1.8in}} 
{\bf\sl Figure 9. 
Here $F^{\circ 2}(v_{10})$ and $F^{\circ 3}(v_{20})$ have return period $2$. 
With the notation above we have $v_{1j_1}=v_{15}$ and $v_{2j_2}=v_{26}$.}
\bigskip
We have proved: 
\bigskip
\bigskip
{\bf 7.13 Proposition.} 
{\it Every subordinated configuration $S$ is tame.} 
\medskip
{\bf Proof.} 
This follows from the discussion in this section.  
\endofproof
\bigskip
\bigskip
 
\centerline{\bf 8. Admissible Projections to Cyclic Sets.} 
\bigskip
\bigskip
Take any maximal {\sl admissible semi-reduced tame} subschema $S'$ of $S$. 
Note that a priori $|S'|$ might be the empty set. 
We will prove that $|S'|=|S|$. 
For this we assume \wlog that if 
$S$ contains a subordinated configuration $S''$; 
then $S'$ has degree at least $3$ (compare Section 7). 
This last hypothesis is needed so we can combine Lemma 8.1 with Lemma 8.7. 
\bigskip
{\bf 8.1 Lemma.} 
{\it Let ${\bf H^*}$ be a (cyclic) \hf 
whose ambient cyclic schema ${\cal M}$ has degree at least 3. 
Let $m \ge 1$ be any positive integer. 
Then there is an extension of ${\bf H^*}$ which contains 
a periodic Julia vertex of return period $m$ which is not post-critical.} 
\medskip
{\bf Proof.} 
(Compare Theorem 1.6) 
The existence of such extension is equivalent to the existence of such orbit in the realization $f \in {\cal P}^{\cal M}$ of ${\bf H^*}$. 
However the only case in which such periodic Julia vertices of 
return period $m$ do not exist, 
is if all of the following conditions are met simultaneously:   
$f \in {\cal P}^{\cal M}$ contains a unique critical point of degree $2$, 
this unique critical point eventually maps to a return $2$ cycle,  
and $m=2$. 
In this particular case 
all these conditions can not be satisfied simultaneously. 
\endofproof
\bigskip
{\bf 8.2 Entrance Points.} 
Let $S_1$ be any consistent subschema of $S$. 
We say that a vertex $v \in |S_1|$ is an {\it entrance point to $S_1$} if 
$v \in F(|S|-|S_1|)$; 
in other words, 
those points in $|S_1|$ which are exposed to the `outside' of $S_1$.
\medskip
Note that because of Lemma 6.5 and Proposition 6.6 all critical cycles contained in $S$ belong to the maximal subschema $S'$. 
Also, by definition no cycle in $S'$ can be superfluous respect to the ambient schema ${\cal M}'$ induced by $S'$. 
For all the other cycles we have the following lemma. 
\bigskip
{\bf 8.3 Lemma.}
{\it Let ${\cal C}$ be a cycle which does not belong to $S'$. 
Suppose $S' + {\cal C}$ is admissible. 
If $S$ is semi-reduced, 
then every entrance point to $S' + {\cal C}$ is saturated respect to ${\cal M}'$.} 
\medskip
{\bf Proof.} 
If this were not the case, 
in the realization of ${\bf H^*}(S')$ we can identify a periodic orbit with ${\cal C}$ (compare Theorem 1.6 and Lemma 8.1). 
As ${\cal C}$ is not superfluous (respect to $S$) 
and $S$ is semi-reduced 
we get a contradiction to the maximality of $S'$ in applying Proposition 6.9. 
\endofproof
\bigskip
{\bf 8.4 Lemma.}  
{\it 
Suppose $v_1,v_2 \in |S|-|S'|$ belong to the same fibre and have distinct images. 
Suppose further that each $v_i$ maps either to an entrance point of $S'$ 
or to a cycle ${\cal C}_i$. 
Under the hypothesis of maximality of $S'$ we have that either 
$v_1$ or $v_2$  is not a critical vertex. } 
\medskip
{\bf Proof.} 
By contradiction, suppose both $v_1$ and $v_2$ are critical. 
We will provide an effective recipe to extend $S'$, 
in contradiction with its maximality. 
\medskip
{\it Case 1: $v_1,v_2$ map to different entrance points.} 
\smallskip 
First note that neither $F(v_1)$ nor $F(v_2)$ can be ends in $S'$. 
In fact, 
if say, $F(v_1)$ is an end, then is not saturated. 
It will follow then by Proposition 6.9 that $S' + \{v_1\}$ is tame in contradiction with maximality.  
Therefore there are vertices say $v_1'$ and $v'_2$ in that same fibre as $v_1$ and $v_2$ such that $F(v_i)=F(v_i')$. 
We will rule out this possibility using the following `folding' technique for \Hf s 
defined as follows. 
\smallskip
Let $[v_1',v_2']$ be the shortest path 
between the vertices $v'_1$ and $v_2'$ in the graph ${\bf H^*}(S')$. 
It follows that the set $F([v_1',v_2'])$ contains the segment 
$[F(v'_1),F(v'_2)]=[F(v_1),F(v_2)]$.\break 
Therefore, after further subdividing the segment $[v_1',v_2']$, 
we may assume there is an edge 
$\ell \subset [v_1',v_2']$ with endpoints say $w_1,w_2$ such that 
$F(\ell)=F([w_1,w_2])=[F(w_1),F(w_2)] \subset [F(v_1),F(v_2)]$. 
To fix ideas, suppose $F(w_1)$ is the closest to $F(v_1)$ among $F(w_1),F(w_2)$ 
(compare Figure 10). 
In this case, in the ordered segment $[w_1,w_2]$ 
we interpolate $v_2$ and $v_1$ in that order, 
making angles of $1/d(v_i)$ between branches.  
\medskip
\centerline{\psfig{figure=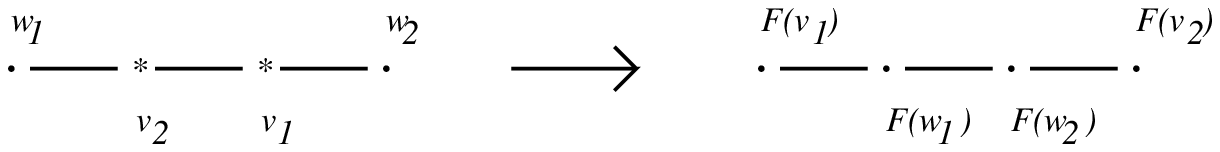,height=1.4in}} 
\noindent
{\bf\sl Figure 10. The folding technique. 
The addition of $v_1$ and $v_2$ does not alter the structure of the image.} 
\bigskip
It follows easily that this new object is a \hf 
which realizes $S' + \{v_1,v_2\}$. 
(We will refer again to this construction in Section 9). 
Now, as we have increased the inner degree of the underlying cyclic schema without adding periodic cycles, 
it follows from Proposition 6.3 that this \hf is tame, 
in contradiction with the maximality of $S'$. 
\medskip
{\it Case 2: $v_1,v_2$ map to disjoint cycles ${\cal C}_1$, ${\cal C}_2$ outside $S'$.} 
\smallskip
In this case the orbits of $v_1$ and $v_2$ form a closed ($\emptyset$-)pseudo-chain. 
As every\break 
($\emptyset$-)closed pseudochain is tame (compare Lemma 6.5 and Section 5), 
we get a contradiction to maximality in applying Proposition 6.6. 
\medskip
{\it Case 3: $v_1$ maps to an entrance point while $v_2$ maps to a cycle ${\cal C}$ outside $S'$.} 
\smallskip
Suppose first that $S'$ has degree at least $3$. 
Then according to Lemma 8.1, 
$S' + {\cal C}$ can be realized and we are in Case $1$. 
\smallskip
Otherwise, if $S'$ has degree $2$, let $v_3$ be the unique critical point in $S'$. 
In particular because $S'$ is semi-reduced, 
$S'$ must be generated by this critical vertex $v_{3}$ and 
therefore must be connected. 
By hypothesis we have then that $v_1$ eventually maps into this component. 
If $v_3$ belongs to a critical cycle, then we are in Case $2$.
Otherwise, if $v_{3}$ is an end, then 
there is a unique pair of points $w_1, w_2$ 
which have the same image. 
It follows from Lemma 3.1 and the fact that $S'$ has degree $2$ that 
$w_1$ and $w_2$ belong to that same fibre as $v_3$; 
furthermore, that same lemma guarantees that they are different from $v_3$. 
Also, only one, say $w_1$, among the pair $w_1,w_2$ belong to that cycle contained in $S'$. 
We distinguish whether $v_3$ and $v_1$ (and therefore $v_2$) belong to the same fibre or not. 
\smallskip
If $v_3$ and $v_1$ belong to the same fibre, 
it follows that 
$S_1= {\cal O}(v_2) + S'$ is a subordinated configuration contained in $S$. 
This is a contradiction with our preliminary hypothesis because $S'$ has only degree $2$.  
\smallskip
If $v_3$ and $v_1$ do not belong to the same fibre, 
we will have in particular that 
$F(v_3) \ne F(v_1)$. 
Also, as $v_3$ is an end, 
$v_3$ is not saturated and therefore Lemma 8.3 shows that it is not an entrance point. 
In particular $F(v_1) \ne v_3$. 
It follows that $F(v_1)$ belongs to the forward orbit of $F^{\circ 2}(v_3)$. 
In this way, there is a unique $v' \ne v_3$ in $|S'|$ such that $F(v')=F(v_1)$. 
In particular $v', v_1,v_2$ all belong to the same fibre. 
Whether $F(v_1)$ belongs to the unique cycle of $S'$ or not, 
the vertex $w_1$ defined above should belong to the forward orbit of $v_1$. 
Thus, arranging the orbits of $v_1,v_2,v_3$ as 
$$
{\bf v_2},F(v_2),\dots, {\bf v_1}, F(v_1),\dots,w_1, {\bf v_3},F(v_2),\dots,v'
$$
we have an ($\emptyset$)-closed pseudo-chain and therefore is tame (compare Section 5).  
Again, because of Proposition 6.6 this is a contradiction with the maximality of $S'$. 
\medskip
{\it Case 4: $v_1,v_2$ map to the same cycle ${\cal C}$ outside $S'$.} 
\smallskip
In this case ${\cal C}$ has return period at least $2$.  
Therefore $\{v_1,v_2\} + {\cal C}$ is tame as shown in Example 3.3.  
The contradiction to maximality follows again from Proposition 6.6. 
\endofproof
\bigskip
{\bf 8.5 Lemma.} 
{\it Under the hypothesis of maximality of $S'$ suppose the set 
${\cal V}=|S|-|S'|$ is non empty. 
Suppose $v \in {\cal V}$ is a critical point. 
Then there exists a critical point $v' \in {\cal V}$ in that same fibre as $v$ 
such that 
$F(v')$ does not belong to $|S'|$.} 
\medskip
{\bf Proof.} 
In fact, 
if ${\cal V}$ contains more than one critical point in that same fibre as $v$, 
this is the content of Lemma 8.4. 
Otherwise $v$ is the unique critical point 
in this fibre which belongs to ${\cal V}$. 
We have to prove that $F(v)$ does not belong to $|S'|$. 
Suppose $F(v)$ does belong to $|S'|$, it follows then from Lemma 8.3 that 
(respect to the ambient schema of $S'$)
$F(v)$ is saturated with respect to the fibre which contains $v$. 
But this is a contradiction with the admissibility of $S$ because 
all critical vertices in that fibre are accounted for 
(these are precisely $v$ and those contained in $|S'|$), 
and $v$ contributes by only $d(v)-1$ to the degree of its fibre but by 
$d(v)$ to the sum over all inverses of $F(v)$ in such fibre. 
This establishes the result. 
\endofproof
\bigskip
{\bf 8.6 Lemma.} 
{\it Under the hypothesis of maximality of $S'$ suppose the set 
${\cal V}=|S|-|S'|$ is non empty. 
Suppose $v \in {\cal V}$ is a critical point. 
If $v$ is the unique critical point in its fibre which belongs to ${\cal V}$, 
then $F(v)$ is not periodic.} 
\medskip
{\bf Proof.} 
By contradiction suppose $F(v)$ belongs to a cycle ${\cal C}$. 
It follows from Lemma 8.5 that ${\cal C}$ and $S'$ are disjoint. 
If $S'+{\cal C}$ is admissible we get a contradiction by applying Lemma 8.5
(more precisely, by copying literaly the proof of Lemma 8.5). 
However, $S'+{\cal C}$ can only fail to be admissible if 
$S'+{\cal O}(v)$ is a subordinated configuration 
(compare Lemma 7.7 and the construction in $\S 7.9$). 
However according to Proposition 7.13, 
this is a contradiction to the maximality of $S'$. 
\endofproof
\medskip
{\bf 8.7 Lemma.} 
{\it Under the hypothesis of maximality of $S'$, 
suppose $v \in {\cal V}=|S|-|S'|$ is a critical point eventually mapping into a cycle ${\cal C}$ outside $S'$. 
Let $w$ be that last point in the orbit of $v$ which is not periodic. 
Then there is a critical vertex $v' \in {\cal V}$ in that same fibre as $w$ 
for which $F(v')$ neither belongs to a cycle nor to $|S'|$.}
\medskip
{\bf Proof.} 
In fact, as $F(w)$ has two preimages, 
by admissibility there must be a critical vertex in that same fibre as $w$. 
If we can find such critical vertex in ${\cal V}$ then the result follows immediately from Lemma 8.5 and Lemma 8.6. 
Otherwise, 
if all critical points in that fibre to which $w$ belongs are assumed to be in $S'$ 
we will derive a contradiction to the maximality of $S'$. 
If $S'$ has degree at least $3$, 
then $F(w) \in {\cal C}$ is not saturated in $S'+{\cal C}$ in contradiction to Lemma 8.3. 
Therefore, $S'$ should have degree 2. 
It follows that the only critical point in $S'$ belongs to that same fibre as $w$ 
and therefore ${\cal O}(v)+S'$ is admissible. 
Furthermore, this subschema represents a subordinated configuration and therefore is tame (compare Section 7). 
This is a contradiction to the maximality of $S'$. 
\endofproof
\bigskip
\bigskip
We can say even more (compare Example 8.11) if we `jump' from fibre to fibre: 
\bigskip
{\bf 8.8 Proposition.} 
{\it Under the hypothesis of maximality of $S'$ suppose the set 
${\cal V}=|S|-|S'|$ is non empty. 
Then ${\cal V}$ contains an ordered sequence of critical vertices $v_{i_1},\dots,v_{i_k}$ belonging to different fibres $W(u_{i_j})$, 
such that 
\smallskip
i) No $F(v_{i_j})$ is an entrance point of $S'$ nor belongs to a cycle; 
\smallskip
ii) That last point $v'_{j}$ in the orbit of $v_{i_j}$ belonging to ${\cal V}-{\cal O}(\{v_{i_1},\dots,v_{i_{j-1}}\})$ or closing that cycle into which $v_{i_j}$ eventually might fall, 
satisfies  
$v'_{j} \in W(u_{i_{j+1}})$ (here $i_{k+1}=i_1$). 
\smallskip
In other words, that portion of the orbits of the $v_{i_j}$ up to $v'_{i_j}$ 
form a closed $|S'|$-pseudo-chain with generator 
$\{v_{i_1},\dots v_{i_k}\}$.} 
\medskip
{\bf Proof.} 
This proposition represents the key step in our development and its proof 
is the most difficult technicality in this work. 
We will find that ordered sequence of critical points in three steps. 
First we will find a distinguish set of critical vertices. 
We will refine this set and the reminding vertices will generate that required sequence. 
\medskip 
{\bf Step 1:} 
As $S$ is semi-reduced and ${\cal V}$ non empty, 
it follows that ${\cal V}$ contains a critical vertex. 
Let $v_{i_1}$ be any critical vertex as guaranteed by Lemma 8.7. 
(This may not be that $v_{i_j}$ guaranteed by the statement of this proposition!) 
Let $u_{i_1}=\phi(v_{i_1})$ and denote by $v'_{i_1}$ 
that last point in the orbit of $v_{i_1}$ which belongs to ${\cal V}$ or 
is {\bf not} periodic if $v_{i_1}$ eventually falls into a cycle outside $S'$. 
To simplify notation denote by 
${\cal L}(v_{i_1})=\{v_{i_1},F(v_{i_1}),\dots,v'_{i_1}\}$ 
that portion of the orbit of $v_{i_1}$ up to $v'_{i_1}$ 
which by Lemma 8.7 contains at least two points. 
\smallskip
We inductively define analogous $v_{i_j}$, $u_{i_j}$, $v'_{i_j}$ 
and ${\cal L}(v_{i_j})$ as follows. 
Suppose these elements have been defined for $j=n$. 
Then we set $u_{i_{n+1}}=\phi(v'_{i_{n}})$. 
If either $u_{i_{n+1}} \in \{u_{i_1},\dots u_{i_{n}} \}$ or 
if $F(v'_{i_{n}}) \in {\cal L}(u_{i_1}) \cup \dots {\cal L}(u_{i_{n-1}})$ 
we stop and continue to step 2. 
\smallskip
Otherwise $F(v'_{i_{n}}) \in |S'|$ or is periodic. 
In either case we take any $v_{i_{n+1}}$ in that same fibre $W(u_{i_{n+1}})$ as $v'_{i_{n}}$ as guaranteed by Lemma 8.7. 
Define $v'_{i_{n+1}}$ as that last element in the orbit of 
$v_{i_{n+1}}$ which belongs to 
${\cal V}-{\cal L}(u_{i_1}) \cup \dots {\cal L}(u_{i_{n}})$
or that last one which is {\bf not} periodic if $v_{i_{n+1}}$ eventually falls into a cycle outside $S'$. 
Finally set 
${\cal L}(v_{i_{n+1}})=\{v_{i_{n+1}},F(v_{i_{n+1}}),\dots,v'_{i_{n+1}}\}$ 
as that portion of the orbit of $v_{i_{n+1}}$ up to $v'_{i_{n+1}}$. 
As the number of fibres is finite, 
it follows that this process must eventually stop. 
\medskip
{\bf Step 2:} 
If $u_{i_{n+1}}=u_{i_1}$ and $F(v_{i_1}) \ne F(v'_{i_n})$, 
then we proceed to Step 3. 
Otherwise, that sequence of critical vertices we are looking for, 
should be found among this set $\{v_{i_1},\dots v_{i_n}\}$. 
However, the election of $v_{i_1}$ as the first element of the sequence seems 
to be inappropriate. 
Therefore we remove it from the list, 
to form a set 
${\cal X}=\{v_{i_2},\dots v_{i_n}\}$. 
However, all elements in the set  ${\cal Y}=\{v'_{i_1},\dots v'_{i_n}\}$ 
convey valuable information. 
For example, the image $F({\cal Y})$ of this set contains by definition 
all those points where the portion of the orbits of two different $v_{i_k},v_{i_l}$
might `collide' outside $S'$. 
Note also that every $v \in {\cal X}$ automatically satisfies condition {\it i)} in the statement. 
Finally also set ${\cal U}=\{u_{i_2},\dots u_{i_n}\}$. 
\smallskip
We start that same process in Step 1 by redefining $v_{i_1}$ to be any vertex in ${\cal X}$. 
In an analogous way we define the elements $v_{i_j}$, $u_{i_j}$, $v'_{i_j}$ and 
${\cal L}(v_{i_j})$. 
However, it is not difficult to check that in this case we can always take this newly defined elements subject to 
$v_{i_j} \in {\cal X}$.  
From this condition it will follow that $u_{i_j}=\phi(v_{i_j}) \in {\cal U}$. 
In fact, it follows easily that if $v_{i_j} \in {\cal X}$ then 
$F(v'_{i_j})=F(v)$ for some $v \in {\cal Y}$. 
Therefore this `new' $v'=v'_{i_j}$ should share the same fibre with an old one. 
The result follows by induction. 
Now by removing one element from ${\cal X}$ at the time if necessary, 
it follows that the sequence must eventually close and we go to Step 3.  
\medskip
{\bf Step 3:} 
If that sequence constructed in Step 2 closes 
(i.e, if $u_{i_1}=u_{i_{n+1}}$ and $F(v_{i_1}) \ne F(v'_{i_n})$),  
the sequence $\{v_{i_1},\dots,v_{i_n}\}$ has that require properties. 
In fact, 
if the orbit of some $v_{i_j}$ eventually falls to a periodic cycle disjoint 
from $S'$ which is avoided by the orbit of those $v_{i_k}$ with $k < j$;  
then, instead of defining $v'_{i_j}$ as 
that last point in the orbit of $v_{i_{j}}$ which is not periodic, 
we define it as the element which closes that cycle 
(in particular $F(v'_{i_{j}})$ denotes the same vertex in both definitions.) 
The result follows. 
\endofproof
\bigskip  
{\bf 8.9 Corollary.} 
{\it That maximal tame subschema $S'$ is non-empty.} 
\medskip
{\bf Proof.} 
In fact, if $S'$ were empty, 
the sequence constructed in Lemma 8.8 generates a tame ($\emptyset$-)closed pseudochain. 
This is a contradiction with the maximality of $S'$. 
\endofproof
\bigskip
{\bf 8.10 Bouquets.} 
Let ${\cal X}=\{v_{i_1},\dots,v_{i_k}\}$ and ${\cal U}=\{u_{i_1},\dots,u_{i_k}\}$ be as in Proposition 8.8. 
We will show that $S''=S'+ {\cal O}({\cal X})$ 
is tame by constructing a modification of ${\bf H^*}(S')$. 
This will be a contradiction with the maximality of $S'$.  
For this we need some preliminary constructions. 
We take for $|S''|-|S'|$ that order 
$$
{\bf v_{i_1}} \prec F(v_{i_1}) \prec \dots \prec v'_{i_1} \prec {\bf v_{i_{2}}} \prec \dots \prec \dots \prec \dots \prec v'_{i_{k},}
$$ 
induced by the definition of pseudo-chain. 
To simplify notation, 
we suppose that $v_{i_1}$ belongs to $W(u_1)$. 
Also, for $u \in {\cal U}$, we denote by $v''_u$ the unique  
$v' \in {\cal Y}$ as in Proposition 8.8 which belongs to the fibre $W(u)$. 
\smallskip
Note that we may assume $F(v_u'') \in |S'|$ for some $u \in {\cal U}$.  
For otherwise, 
we will have a closed ($\emptyset$-)pseudochain, 
which according to Lemma 6.5 is a contradiction to the maximality of $S'$. 
Thus, as the pseudo-chain is closed, rearranging the indices if necessary 
we assume that $F(v_1'') \in |S'|$. 
We proceed to construct an angled tree $B_u$ for every $u \in |{\cal M}|$. 
(This graph would look like more as a bouquet than as a forest.) 
Later this `flowers' will be grafted along the tame cycle of ${\bf H^*}(S')$  
to construct a Hubbard Forest ${\bf H^*}(S'')$. 
\medskip
The construction of $B_1$ is standard. 
Let $w_1=v_{i_1},\dots ,w_m$ be the ordered collection of points in $W(u_{1}) \cap (|S''|-|S'|)$. 
Note that this $m$ is always bigger than one. 
We join all $w_\alpha$ to a vertex $q_1$ 
following that induced order in the pseudo-chain 
and making angles of $1/m$ between consecutive edges 
(a usual star configuration). 
Finally we join a new vertex $p_1$ to $w_1$ making an angle of $1/d(w_1)$ with this star configuration. 
(Note that by definition $d(w_1)=d(v_{i_1}) > 1$.) 
This `flower' is $B_1$. 
(Below, this point $p_1$ would be identified with the {\it distinguish tame point $p_1$} of ${\bf H^*}(S')$.)
\smallskip
Next note that by definition $w_m=v''_1$. 
Thus, $w_m$ would `map' into the forest ${\bf H^*}(S')$. 
As the star itself has no critical vertices other than at the ends, 
the image of the star should also be a star in order to satisfy the angle condition. 
The ends of this second star (with center $q_2$) include 
$F(w_1),\dots,F(w_{m-1})$ in that order. 
At the missing end (corresponding to $F(w_m)$) we include a point $p_2$. 
($F(w_m)$ should not be in this configuration. In fact, 
by assumption it is an entrance point to $S'$, 
and therefore has already a representative in ${\bf H^*}(S')$.) 
Now if $u_2=F_{\cal M}(u_1) \in {\cal U}$, 
then necessarily $v''_{u_2}=F(w_j)$ for some $j=1,\dots m-1$. 
In this case we interpolate the unique critical vertex $v \in {\cal X} \cap W(u_2)$ 
between the center $q_2$ and $v''_{u_2}=F(w_j)$ making angles of $1/d(v)$ 
between branches (compare Figure 11). 
Note that by definition of closed pseudo-chain we necessarily have 
$F(v) \ne F(v''_{u_2})$. 
Furthermore, the inclusion of $v$ in this place is done so that we can comply with several conditions. 
For example $v$ will take the place of $v''_{u_2}$ in the construction of the next star. 
In fact, `the chain' is broken at this point $v''_{u_2}$, 
so that $F(v''_{u_2})$ 
will not follow that `order of the star', 
$F(v''_{u_2})$ might even be in ${\bf H^*}(S')$. 
Furthermore, this $F(v''_{u_2})$ usually has at least two preimages and therefore 
a critical point should be interpolated between (compare Lemma 3.2). 
Otherwise if $u_2 \not\in {\cal U}$, 
there is nothing to worry about. 
In either case the resulting configuration is defined to be $B_2$. 
\smallskip
Now we continue inductively this process  for all $u_{i}$. 
In other words, we push the star to the next component. 
We artificially complete the missing edge by joining the center of the star to a point $p_u$. 
Finally whenever there is vertex at which the order of the chain is broken, 
we interpolate between this vertex and the center that unique critical vertex 
$v \in {\cal X} \cap W(u)$. 
Note that the stars `rotates' only when mapping from the first to the second component. 
This guarantees the compatibility of the dynamics in the last branch $B_{r}$. 
Also note that at every vertex where more than one edge meets, 
by construction the angle condition is satisfied. 
\medskip
Now we modify the Forest ${\bf H^*}(S')$ to allow place for the `flowers' $B_u$ as follows. 
Because ${\bf H^*}(S')$ is tame, 
every component of ${\bf H^*}(S')$ has a distinguish periodic point $p_u$ (of return period 1 and rotation number zero) 
at which say $n \ge 1$ edges meets. 
The angle between consecutive edges is there $1/n$. 
We modify this angle to $1/(n+1)$ 
and allow space for an extra branch where we `glue' $B_u$. 
Now, there might be points in ${\bf H^*}(S')$ which eventually map to $p_u$. 
But because ${\bf H^*}(S')$ is tame, 
it follows that there is no need to modify the angle between edges at such inverses. 
With the induced dynamics and degree from $S$ at new vertices 
(and by defining $q_{u_\alpha}$ to be of degree one with $f(q_{u_\alpha})=q_{u_{\alpha+1}}$), 
this abstract Hubbard Forest realizes $S''$. 
Furthermore, it is easy to see that 
$p_0 \mapsto \dots \mapsto p_r=p_0$ still defines a tame cycle. 
In fact, no vertex $v \in |S''|$ other than those in $S'$ map to the cycle. 
But at these vertices we know condition (T) is satisfied.  
We have then that $S''$ is tame, in contradiction to the maximality of $S'$.  
\medskip
{\bf 8.11 Example.} 
Consider the admissible schema $S$, with diagram 
\smallskip
\centerline{\psfig{figure=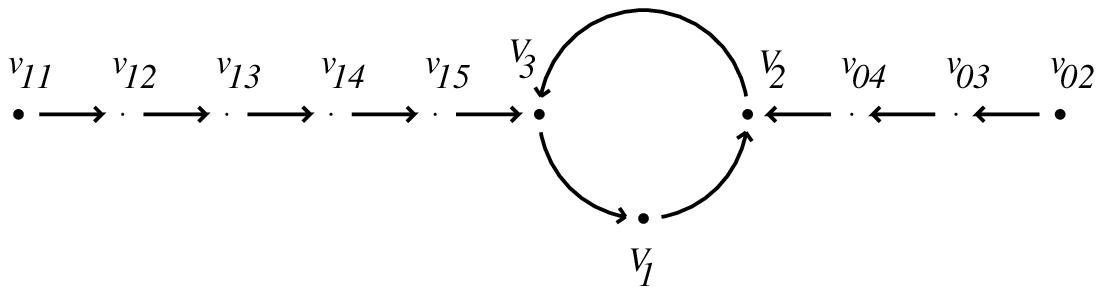,height=1.5in}} 
\noindent
and fibres 
$$
\eqalign{
W(1) &=\{V_1,v_{04},v_{11},v_{14}\}\cr
W(2) &=\{V_2,v_{02},v_{12},v_{15}\}\cr
W(3) &=\{V_3,v_{03},v_{13}\}. 
}
$$
\smallskip
The cycle ${\cal C}: \hskip 0.05in V_1 \mapsto V_2 \mapsto V_3 \mapsto V_1$ is easy to realize as was shown in Section 4. 
When this is done, 
we have that ${\cal C}$ is saturated, 
so we need to proceed as in $\S$8.10. 
Here $v_{11},v_{02}$ are vertices as guaranteed by Proposition 8.8. 
Thus, we must take for $S-{\cal C}$ the order 
$$
{\bf v_{11}} \prec v_{12} \prec v_{13} \prec v_{14} \prec v_{15} \prec {\bf v_{02}} \prec v_{03} \prec v_{04}, 
$$
to form a sequence $v_{11},v_{14},v_{04}$ in $W(1) \cap (S-{\cal C})$. 
According to $\S$8.10 we should then construct a `bouquet' by starting with this order. 
(Compare Figure 11). 
\smallskip
\centerline{\psfig{figure=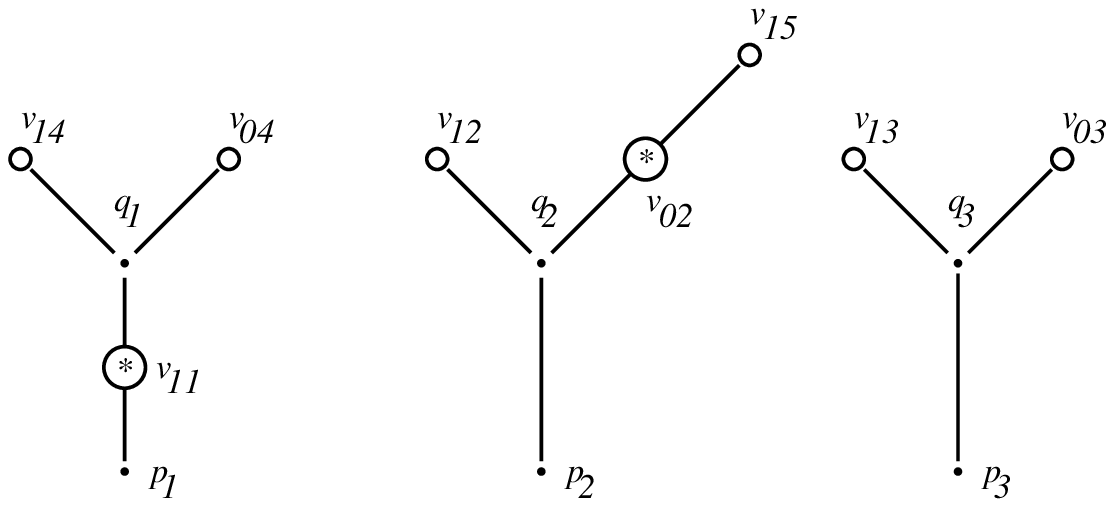,height=1.9in}} 
\noindent
{\bf\sl Figure 11. A bouquet. 
Note that $v_{04}$ and $v_{15}$ map outside the bouquet. 
Also the dynamics at $v_{11}$ `rotates' clock-wise the star $1/3$ turns. 
These two facts guarantee the angle condition is satisfied.} 
\smallskip
When this is done, 
we glue this bouquet around $p_1$, $p_2$, $p_3$; 
and this \hf realizes $S$. 
(Compare Figure 12). 
\smallskip
\centerline{\psfig{figure=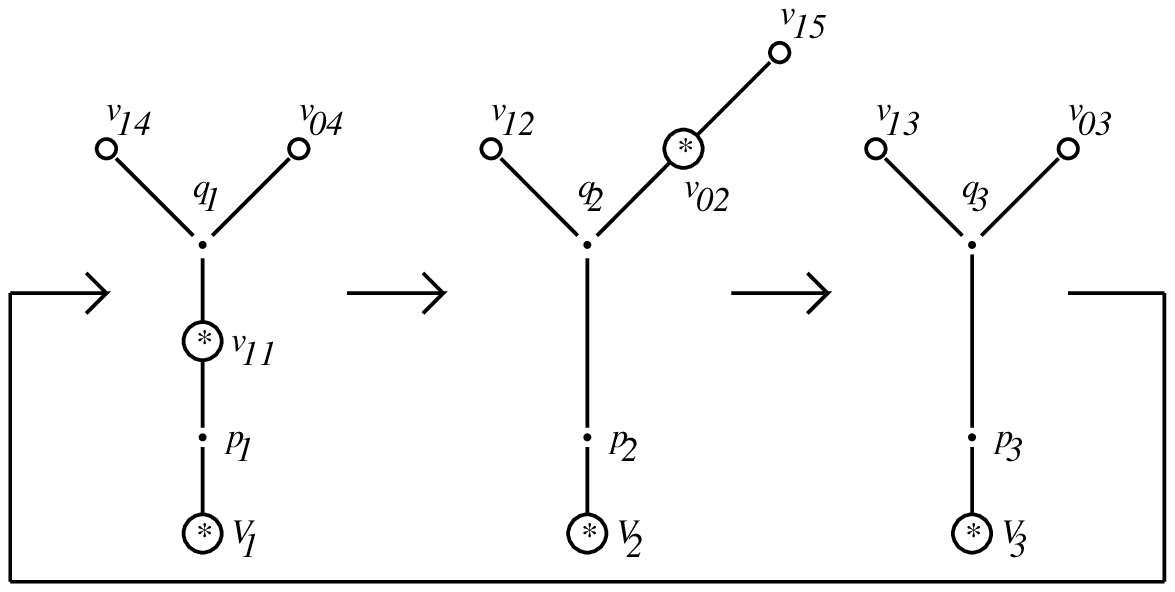,height=2.4in}} 
\centerline{\bf\sl Figure 12. This forest realizes the schema $S$.}
\bigskip
\bigskip
Putting together the results of this section we have the following result. 
\bigskip
{\bf 8.12 Theorem.} 
{\it Any admissible semi-reduced schema $S$ with underlying cyclic\break
ambient schema is tame.} 
\endofproof
\bigskip
\bigskip

\centerline{\bf 9. Coda.} 
\bigskip
\bigskip
The proof of Theorem B is now routine. 
Let 
${\cal C}: \hskip 0.15in u_0 \mapsto u_1 \mapsto \dots \mapsto u_r=u_0$ 
be that unique cycle in ${\cal M}$. 
We realize first that part of $S$ which  projects into this ${\cal C}$. 
In other words we restrict $S$ to the set $\cup_{\alpha=0}^{r-1} W(u_\alpha)$ 
to form a subschema $S_1$. 
To this schema $S_1$ we can associate a semi-reduced schema $S_2$ by inductively removing all superfluous cycles and ends which are not critical. 
By Theorem 8.12 this semi-reduced schema $S_2$ 
can be realized as a Hubbard Forest ${\bf H^*}(S_2)$. 
As the subschema $S_1$ is admissible, 
it follows that the polynomials that realize $S_2$ would also realize $S_1$ 
(with an extended set of `marked point'). 
In other words an extension ${\bf H^*}(S_1)$ of ${\bf H^*}(S_2)$ realizes $S_1$ as a Hubbard Forest. 
\smallskip
Next we inductively add an extra vertex $u \in |{\cal M}|$ at the time 
to that part of the ambient schema ${\cal M}$ which is already realized until it is complete. 
For this, we take $u \in |{\cal M}|$ not belonging to the cyclic part of ${\cal M}$. 
We suppose that $H_{F_{\cal M}(u)}$ has been constructed and proceed to define $H_u$. 
We start by defining a graph $H_u$ as a copy of $H_{F_{\cal M}(u)}$. 
In this angled tree, we consider a maximal set of vertices in the fibre $W(u)$ which can be embedded as vertices. 
We ask at this step only for compatibility with the dynamics, 
in other words, the degree at each vertex is not copied. 
The degree at this points is that induced from $S$. 
If among this maximal set there are critical points, 
we proceed to modify the angle function of the tree so that the angle condition is satisfied at those points. 
We extend this new angled tree to be a homogeneous covering of $H_{F_{\cal M}(u)}$ 
(compare Lemma A.4 in the Appendix), 
and look again for a maximal set of extra vertices that can be embedded and so on. 
It may happen that at one stage of this construction no new vertices in $W(u)$ can be embedded. 
If this is the case and there are still vertices in $W(u)$ not embedded, 
the condition of admissibility shows there are at least $2$ critical points with different image, 
which are yet to be embedded. 
In this case we appeal to the folding technique for coverings 
to include this two critical vertices simultaneously 
(compare Lemma 8.4, Case 1). 
Thus, after a finite number of steps we are done. 
This completes the proof of Theorem B. 
\bigskip
\bigskip

\centerline{\bf Appendix A: Hubbard Forests.} 
\bigskip
\bigskip
{\bf A.1. Angled Trees.} 
By an {\it (angled) tree H=(T,V,\angle)} will be meant a finite connected acyclic $m$-dimensional simplicial complex where $m=0$ or $m=1$, 
together with a function 
$\ell,\ell' \mapsto \angle(\ell,\ell')=\angle_v(\ell,\ell') \in {\bf Q/Z}$ 
which assigns a rational modulo 1 to each pair of edges $\ell,\ell'$ which meet at a common vertex $v$. 
This angle $\angle(\ell,\ell')$ should be skew-symmetric, with $\angle(\ell,\ell')=0$ if and only if $\ell=\ell'$, 
and with $\angle_v(\ell,\ell'')=\angle_v(\ell,\ell')+\angle_v(\ell',\ell'')$ whenever three edges are incident at a vertex $v$. 
Here $V$ represents the set of vertices and $T$ the underlying topological tree 
(or a point if $m=0$). 
Such an angle function determines a preferred isotopy class of embeddings of $T$ into ${\bf C}$, 
in which the induced cyclic order of edges at a vertex $v$ is preserved. 
\bigskip
{\bf A.2. Definition.} 
Let $H_i=(T_i,V_i,\angle_i)$ $(i=1,2)$ be angled trees. 
By an\break
{\it abstract $H_1$-covering of $H_2$} will be meant a function 
$f:V_1 \to V_2$, together with 
a local degree 
$\delta_f:V_1 \to {\bf Z}$ which assigns a positive integer $\delta_f(v)\ge 1$ to every vertex $v \in V_1$. 
These functions are subject to the following restrictions. 
\smallskip
Whenever $v$ and $v'$ are the endpoints of a common edge $\ell$ in $H_1$, 
we require that\break
$f(v) \ne f(v')$. 
Next, extend $f$ to a map $f:H_1 \to H_2$ which carries each edge homeomorphically onto the shortest path joining the images of its endpoints. 
We require then that 
$\angle_{f(v)}(f(\ell),f(\ell'))=\delta_f(v)\angle_v(\ell,\ell')$. 
(In this case $f(\ell)$ and $f(\ell')$ are incident at the vertex $f(v)$ where the angle is measured.) 
\smallskip
By definition a vertex {\it $v \in V_1$ is critical} if $\delta(v)>1$ and 
{\it non-critical} otherwise. 
The number $deg(f)=1+\sum_{v \in V_1}(\delta(v)-1)$ is the {\it degree of $f$}. 
\medskip
A particular case is when $deg(f)=1$. 
In this case we say that $H_2$ is an {\it extension of $H_1$}. 
Here we will think of $H_1$ as a subset of $H_2$ 
(compare lemma A.4 where this abuse of notation is used). 
\bigskip
{\bf A.3. Homogeneous Coverings.} 
We say that an abstract $H_1$-covering of $H_2$ is {\it homogeneous} if for every $v \in V_2$ we have
$$
deg(f)=\sum_{v' \in V_1:f(v')=v} \delta_f(v'). 
$$
In other words every vertex in $V_2$ has as many inverses counting `multiplicity' by $\delta_f$ as the degree of $f$. 
\medskip
In general it is easier to work with homogeneous coverings. 
In this context this represents no loss of generality as the following lemma shows. 
\bigskip
{\bf A.4. Lemma.} 
{\it Let $(f,\delta_f)$ be an abstract 
$H_1$-covering of $H_2$ of degree $n>0$. 
Then there exists an extension angled tree $H_1'$ of $H_1$ and 
$(f',\delta_{f'})$ a homogeneous  
abstract $H'_1$-covering of $H_2$ also of degree n, 
such that 
$f(v)=f'(v)$ and $\delta_f(v)=\delta_{f'}(v)$ for every $v \in V_1$. 
Furthermore such extension is unique.
} 
\bigskip
{\bf A.5 Corollary.} 
{\it Let $(f,\delta_f)$ be an abstract $H_1$-covering of $H_2$ of degree $n>0$. 
Then every point $v \in V_2$ has at most $n$ inverses counting multiplicity.} 
\bigskip
\bigskip
{\bf A.6. Definition.} By an {\it abstract Hubbard Forest ${\bf H^*}$ 
with underlying reduced mapping schema ${\cal M}=(|{\cal M}|,F,w)$}
will be meant 
a collection $H_u$ ($u \in |{\cal M}|$) of angled trees and 
a collection $(f_u,\delta_u)$ of abstract $H_u$-coverings of $H_{F(u)}$ of degree $d(u)=w(u)+1$; 
satisfying the following (usual) conditions. 
\smallskip
Let $V=\bigcup_{u \in |{\cal M}|} V_u$. 
We say that vertex $v \in V$ is a {\it Julia vertex} if  
there are no periodic critical vertices in its forward orbit. 
We have a well defined function $f:V \to V$ between vertices 
as defined by $f(v)=f_u(v)$ whenever $v \in V_u$. 
We require then that 
for any periodic Julia vertex $v$ where $m$ edges meet, 
the angle between two edges incident at $v$ 
should be a non trivial multiple of $1/m$. 
(A vertex which is not Julia will be called a {\it Fatou vertex}.)
\smallskip
Let $v,v'$ be vertices in the same tree $H_u$, 
we define the {\it distance $d_{\bf H^*}(v,v')$} between these vertices 
as the number of edges between $v$ and $v'$. 
We require then that 
for any pair of different periodic Julia vertices $v$ and $v'$ 
belonging to the same tree $H_u$ there is a $k \ge 0$ such that 
$d_{\bf H^*}(f^{\circ k}(v),f^{\circ k}(v'))>1$. 
\medskip
{\bf Remark.} 
If $f \in {\cal P}^{\cal M}$ is post-critically finite, 
then any invariant set $S$ containing the critical set $\Omega(f)$ naturally defines an abstract Hubbard Forest ${\bf H^*}_{f,S}$. 
\bigskip
\bigskip
{\bf A.7. Theorem.} 
{\it For every abstract Hubbard Forest ${\bf H^*}$ 
with underlying ambient schema ${\cal M}$, 
there is a Post-critically Finite $f \in {\cal P}^{\cal M}$ 
and an invariant set $S$ containing the critical set $\Omega(f)$ which realizes ${\bf H^*}$, 
i.e, ${\bf H^*} \cong {\bf H^*}_{f,S}$. 
Furthermore, 
$f$ is unique up to component-wise affine change of coordinates.} 
\bigskip
\bigskip


\input psfig

\newif\ifboxfigure      
\boxfigurefalse

\def\BoxIt#1#2{
	\vbox{\hrule
	\hbox{\vrule\kern#2\vbox{\kern#2#1\kern#2}\kern#2\vrule}
		   \hrule}}

\def\insertRaster #1 pixels #2  by #3 scaled #4 {
			 \hbox {%

			 \hss
			 \RasterBox {#1} {#2} {#3} {#4}
			 \hss
			 }%
}

\def\RasterBox #1 #2 #3 #4{


\dimen0=#2pt   
\dimen1=#3pt
\divide\dimen0 by 11076 
\multiply\dimen0 by 10 
\divide\dimen1 by 11076
\multiply\dimen1 by 10
 \dimen2=#3pt
 \divide\dimen2 by 933
 \dimen3=#2pt
 \divide\dimen3 by 933
\setbox4=\hbox to #4\dimen0{
 \vbox to #4\dimen1{
 \vss
 \psfig{figure=#1,height=#4\dimen2,width=#4\dimen3}
 }
 \hss
 }
 \ifboxfigure\BoxIt{\box4}{0pt}
 \else\box4
 \fi
 }


\bigskip
\bigskip
\centerline{\bf Appendix B: Ambient Schemata versus Regular Schemata.} 
\bigskip
\bigskip
In this appendix we will show how mapping schema naturally appears 
in the context of dynamics of polynomials. 
We start with the mapping schema  $(2,1)$ which by definition contains two  points, one of weight 2 and the other of weight 1; 
mapping to each other. 
Depending in the context, 
this mapping schema can represent several things: 
\medskip
{\bf 1. As Ambient Schema.} As ambient schema it represents two copies of the complex plane, 
``mapping to each other" with $2$ and $1$ critical points respectively. 
In other words, we think, of this schema, as the proper homotopy class of 
any pair of polynomials $(P_1,P_2)$, so that 
$P_1:{\bf C}_1 \to {\bf C}_2$ has degree $3$ ($=2+1$) and 
$P_2:{\bf C}_2 \to {\bf C}_1$ has degree $2$ ($=1+1$).  
As we are only interested in dynamical behaviour, 
it is not difficult to show that we may restrict ourselves to centered monic 
polynomials $P_1$ and $P_2$ (compare [M]). 
\medskip
{\bf 2. As the restriction to the postcritical set of the dynamics of a single polynomial.} 
For example the polynomial 
$$
P(z)=z^4-6b^2z^2-{8 / 9}z-{5b / 3} \hskip 1in 
\hbox{where $b^3=-{1 / 9}$ } 
$$ 
has critical points $b,b,-2b$ (i.e, $z_0=b$ is a double critical point) 
which satisfy 
$P(b)=-2b$ and $P(-2b)=b$. 
In this way, the restrictriction of the dynamics to the postcritical set 
can be described also as the cyclic schema $(2,1)$. 
\bigskip
\centerline{
\insertRaster 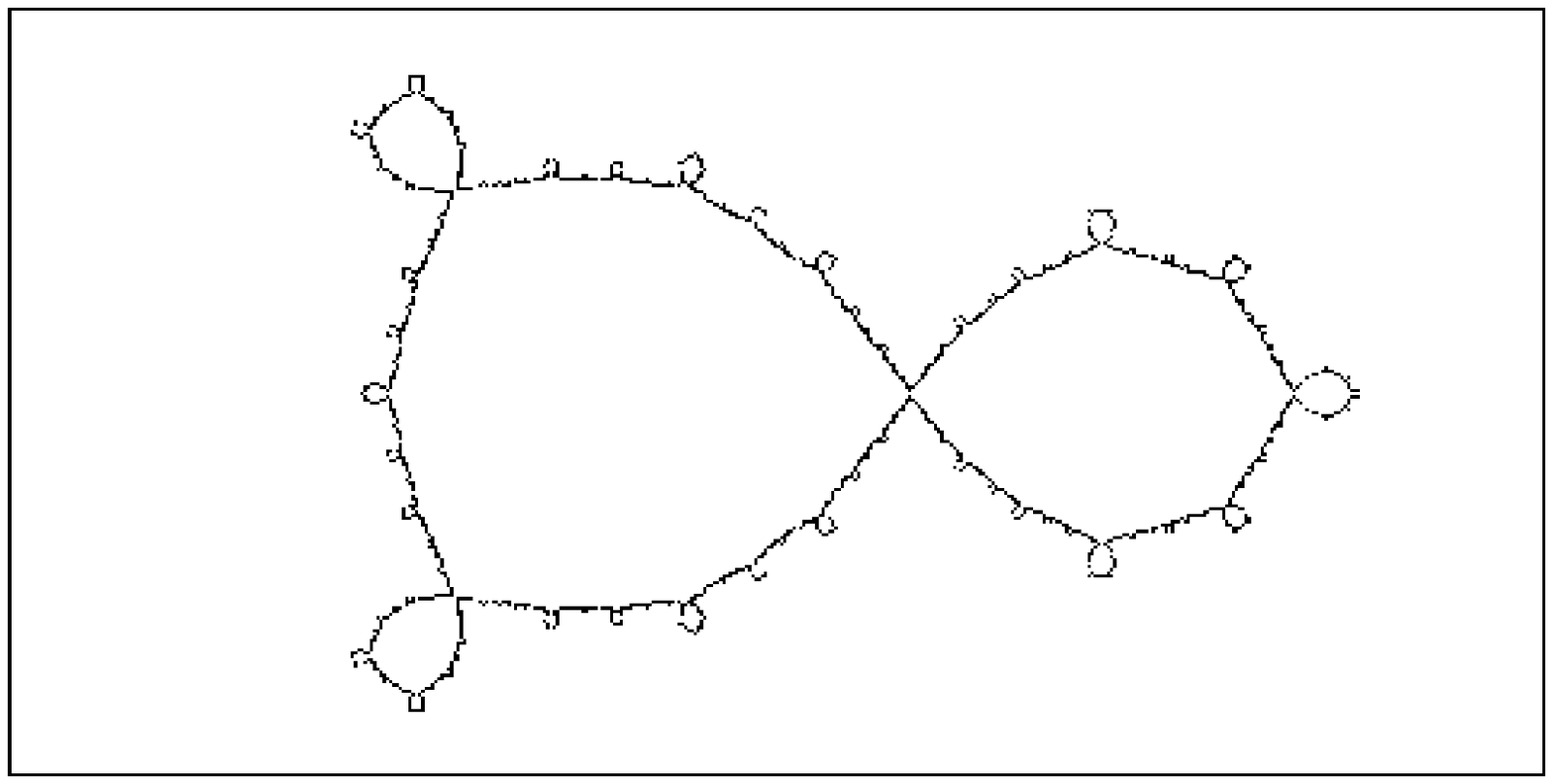 pixels 480 by 240 scaled 550 } 
\smallskip
\centerline{\sl A polynomial with cyclic post-critical schema $(2,1)$}
\bigskip
{\bf 3. As the reduced schema of the dynamics of the postcritical set. }
Most times at some points in the orbit of the critical set the polynomial is 1-1. 
In this case, 
for several purposes the dynamics at these points is simply ignored. 
For example the postcritically finite polynomial with Hubbard Tree
\bigskip
\centerline{\psfig{figure=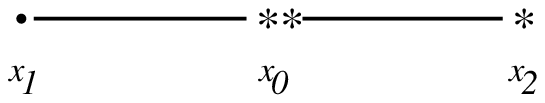,height=1.0in}}
\noindent
has associated cyclic schema ${\cal S}=(2,0,1)$. 
However, if we disregard the dynamics at non critical points, 
we will obtain the associated reduced schema $\bar {\cal S}=(2,1)$. 
\bigskip
\bigskip
{\bf 4.} Let us work a completely different example. 
We consider the cyclic schema ${\cal M}=(1,1)$ containing two critical points
each of weight 1, and which map to each other. 
In practice this schema (as any other will under the appropiate conditions) 
appears in three ways. 
First as an ambient schema; i.e, this describes the space of polynomial maps 
${\cal P}^{(1,1)}$ as defined in section 1.1.2 
(see also $\S$1 in this appendix). 
Then as a postcritically finite polynomial of degree 3 whose two critical points are interchanged by iteration. This polynomial is $f(z)=z^3-1.5z$. 
\bigskip
\centerline{
\insertRaster 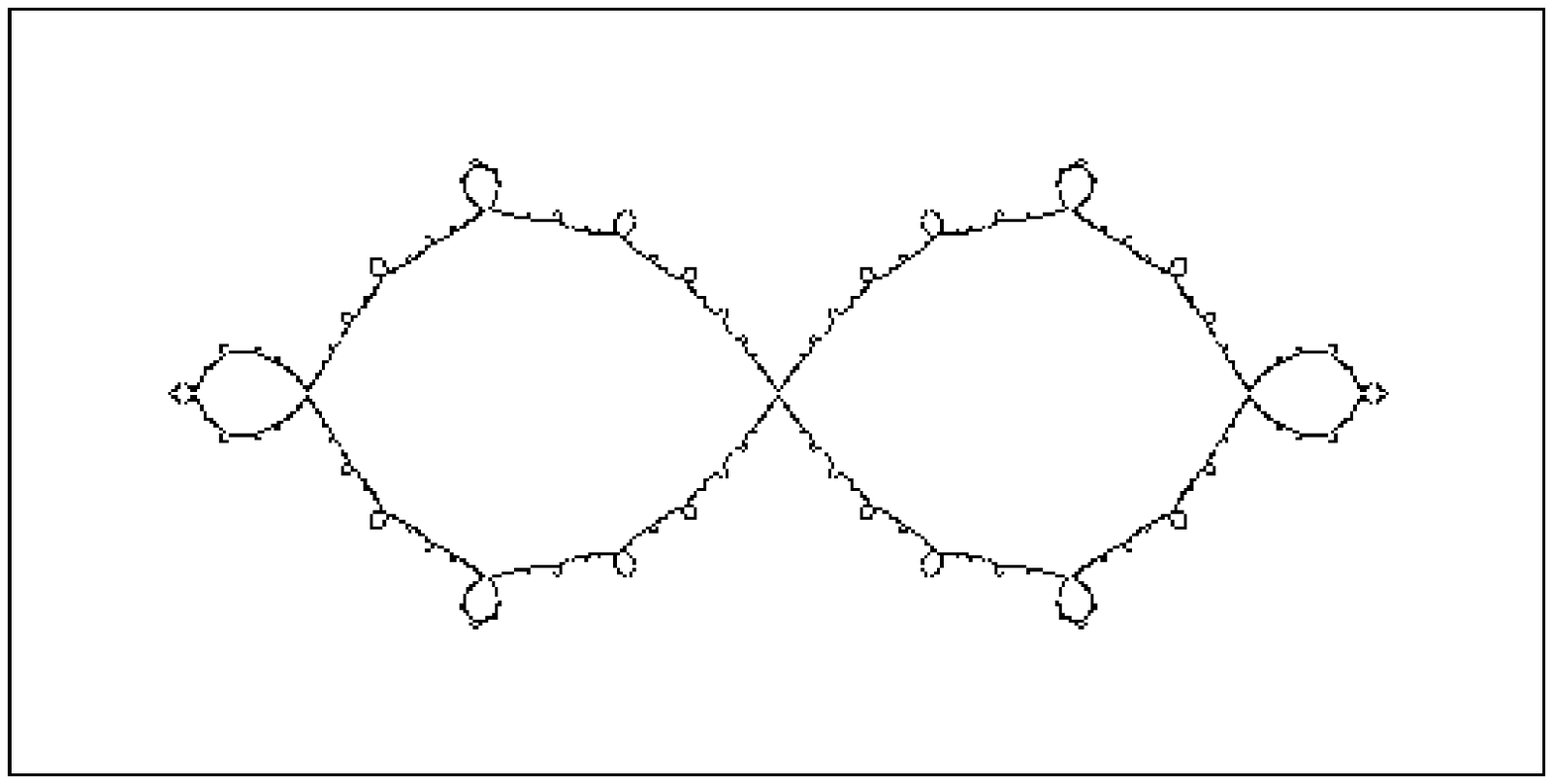 pixels 480 by 240 scaled 550 }  
\medskip
\centerline{\sl  Critical points interchanged by iteration.} 
\medskip 
\noindent
And in the third place as the schema associated to the critical orbits 
of the polynomial map $(z^2,z^2) \in {\cal P}^{(1,1)}$
(this polynomial function maps the first copy of the complex numbers to the second 
by $z^2$, and the second to the first also by $z^2$; 
Therefoe it is the ``preferred map" in the space ${\cal P}^{(1,1)}$). 
There are several things that should  be said about this map. 
First it should be clear that the restriction of the dynamics to the critical set is described by ${\cal M}$. 
However, this map itself lives in the space ${\cal P}^{\cal M}$ and there 
should be a distinction when using ${\cal M}$ as the schema of the postcritical set or as the ambient schema. 
(Evidently this is a very particalar case: in general a postcritically finite 
map $f \in {\cal P}^{\cal M}$ will have associated schema different from ${\cal M}$ as we will see below.) 
In second place, the Julia set of this particular map is the disjoint union 
of two circles. 
\medskip
The relation between the ambient schema ${\cal M}$ 
(or more preciselly between the space ${\cal P}^{\cal M}$) 
and the degree three 
polynomial with associated schema ${\cal M}$
can be heuristically described as follows. 
An ambient schema ${\cal M}$ consists of two copies of the complex numbers,
while our polynomial map has two distinct Fatou componets each one associated with a corresponding critical point. 
We think of each of these Fatou components as canonically corresponding 
to the respective copy in the ambient space.  
Therefore as the parameters change in ${\cal P}^{\cal M}$ 
we shall expect certain similarities in the change of parameters of the degree 
three polynomial 
(examples of this will follow shortly). 
Technically our hope is the following. 
By ``renormalizing" each Fatou component independently 
(but following the $(1,1)$  pattern) 
it is not difficult to check that the obtained renormalization will correspond 
to the map $(z^2,z^2) \in {\cal P}^{\cal M}$ described above. 
In this way we will have that every ``renormalization with the $(1,1)$ 
pattern will correspond to a map in ${\cal P}^{\cal M}$. 
\bigskip
Now we can ask how to classify the possible postcritically finite maps in the space 
${\cal P}^{\cal M}$. 
This was done the author in [P3] and was the main reason for the introduction of Hubbard Forests. 
A second, more general question is which associated schemata ${\cal S}$ 
can correspond to a postcritically finite map in ${\cal P}^{\cal M}$. 
In fact, the purpose of this paper is to give an answer to this question. 
\bigskip
We conclude this appendix with several examples showing the similarities 
between certain maps in ${\cal P}^{\cal M}$ and degree three polynomials close to our $(1,1)$ map. 
\smallskip
The map $(z^2-1,z^2-1) \in {\cal P}^{\cal M}$ has associated schema 
$(1,0)+(1,0)$
and its Julia set consists of two copies of the well known Julia set of $z^2-1$. 
The polynomial $f(z)=z^3-az$ where $a \approx 2.12132$
also has associated mapping schema $(1,0)+(1,0)$ and is close to the 
$(1,1)$ map.  
\medskip
\centerline{
\insertRaster 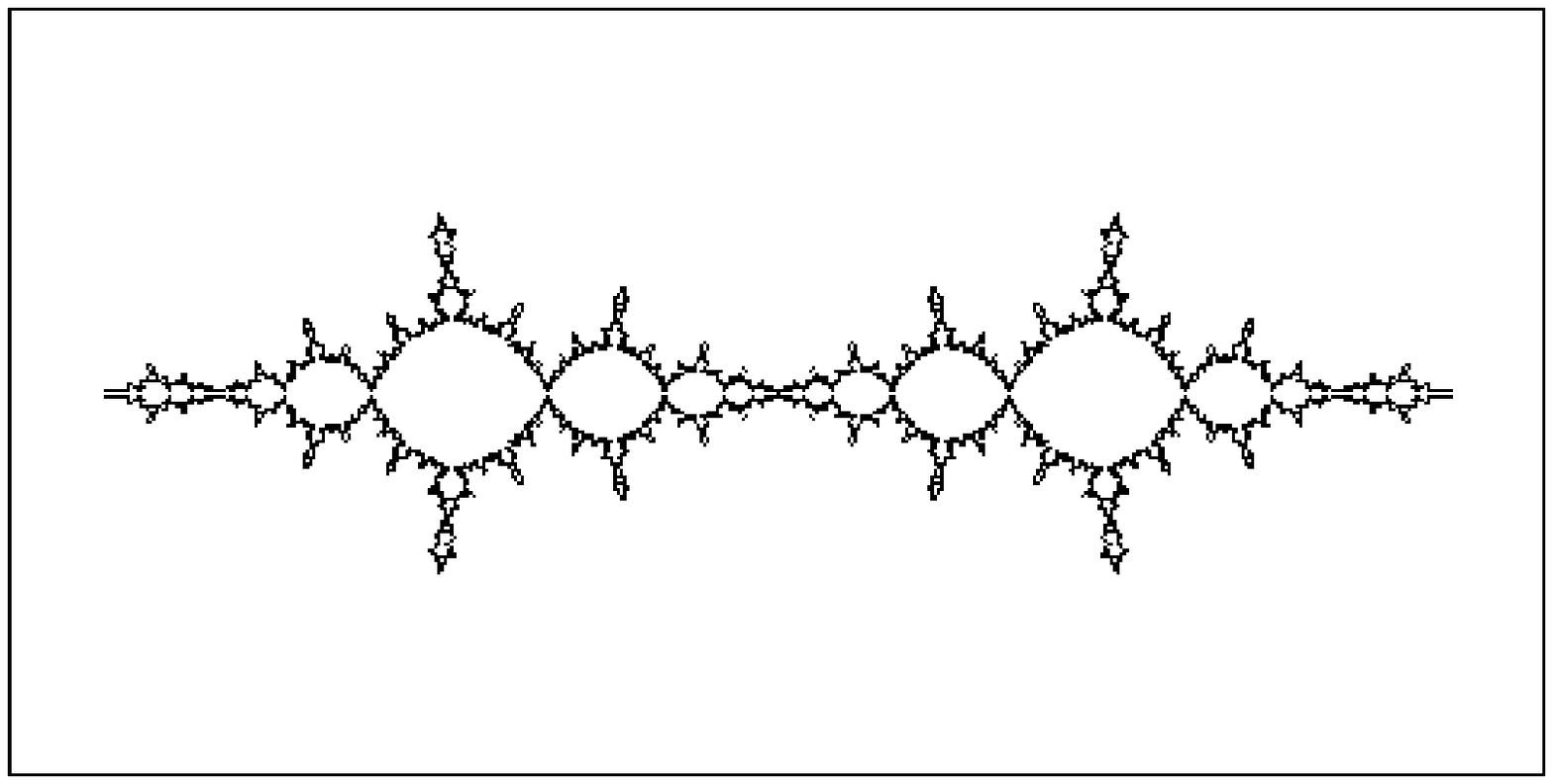 pixels 480 by 240 scaled 550 } 
\medskip 
\centerline{\sl Two copies of the ``Basilica" with extra decorations.} 
\bigskip
Now consider the schema $(\{\{\{1\}1\}0\}0,0)$ 
(i.e, there are two critical points, the first map to the second. 
This second critical point maps in two iterations to a period two orbit. 
The postcritically finite map 
with associated Hubbard Forest 
\bigskip
\centerline{\psfig{figure=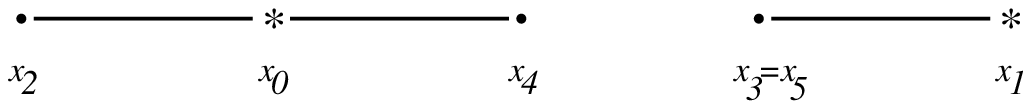,height=1.0in}}
\medskip
\noindent 
corresponds to this last schema. 
In fact the map $(z^2,z^2+c) \in {\cal P}^{\cal M}$ (where $c^3=-2$) 
realizes this Hubbard Forest. 
On the other side, the polynomial $f(z)=z^3-az+b$ 
(where $a \approx 2.004865$ and $b \approx -0.275148$) 
also realizes this last schema and is `close' to the required polynomial. 
\bigskip
\centerline{
\insertRaster 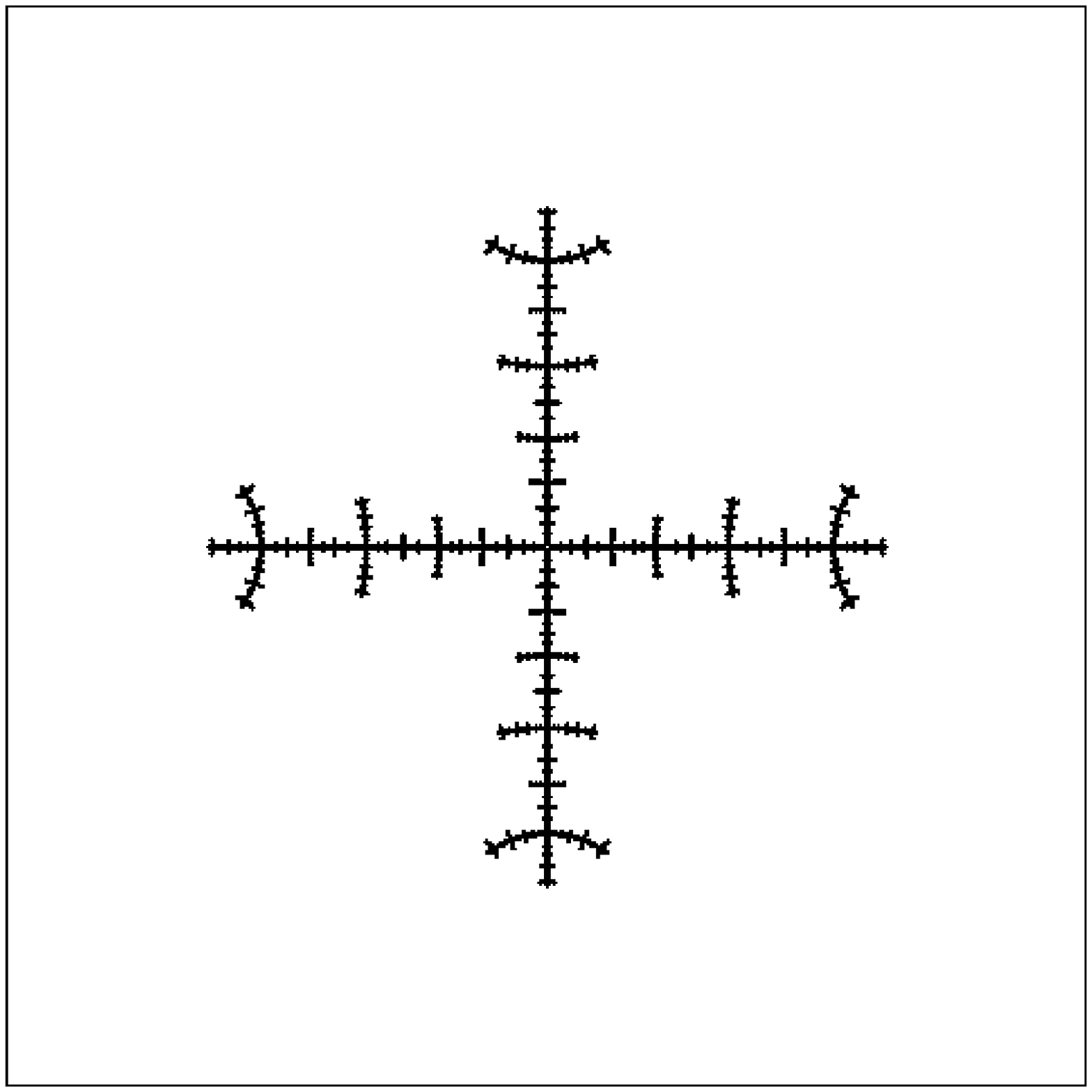 pixels 480 by 480 scaled 375 
\hskip 0.5in 
\insertRaster 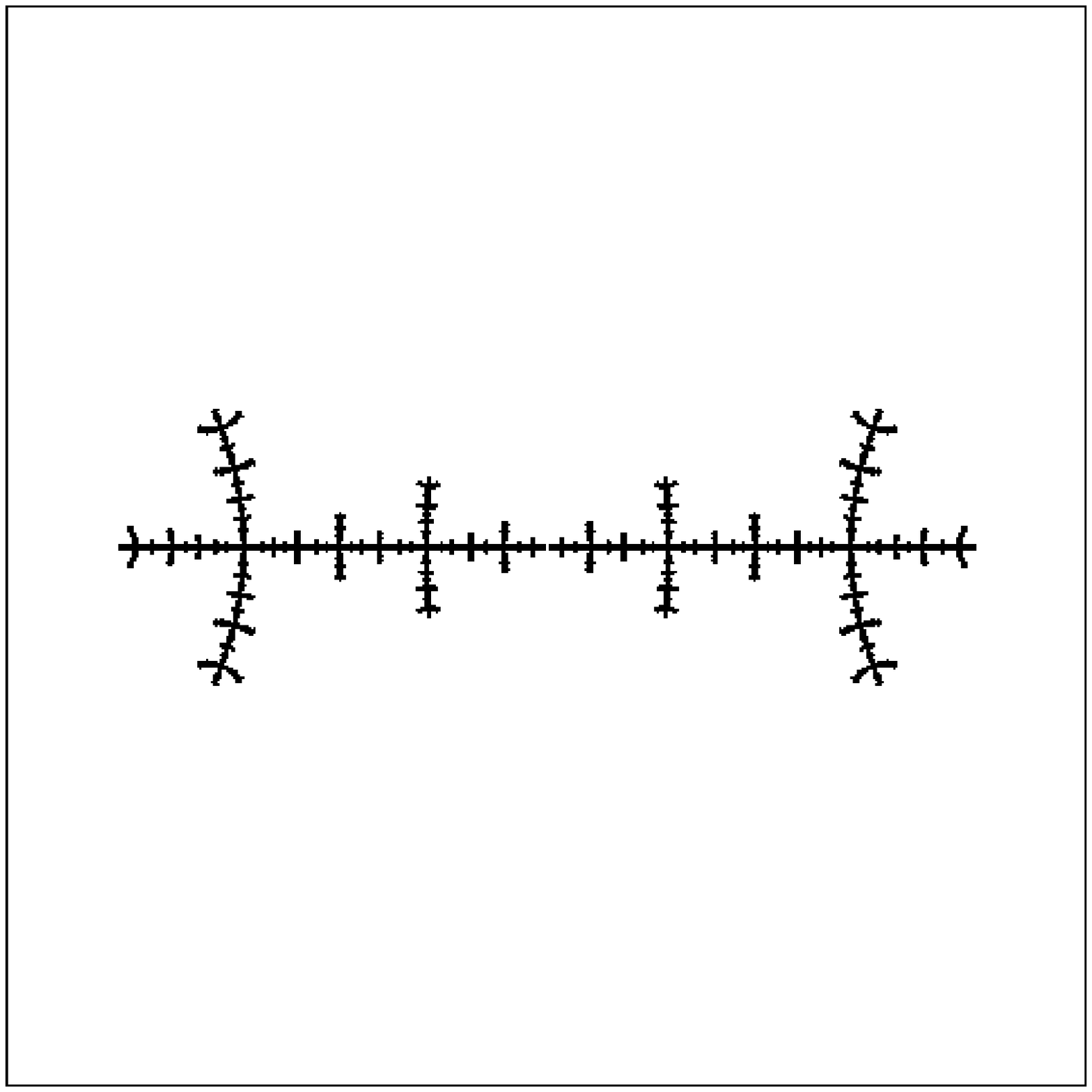 pixels 480 by 480 scaled 375 }
\smallskip
\centerline{\sl A pair of Polynomial maps with ambient schema $(1,1)$.}
\bigskip
\bigskip
\centerline{ 
\insertRaster 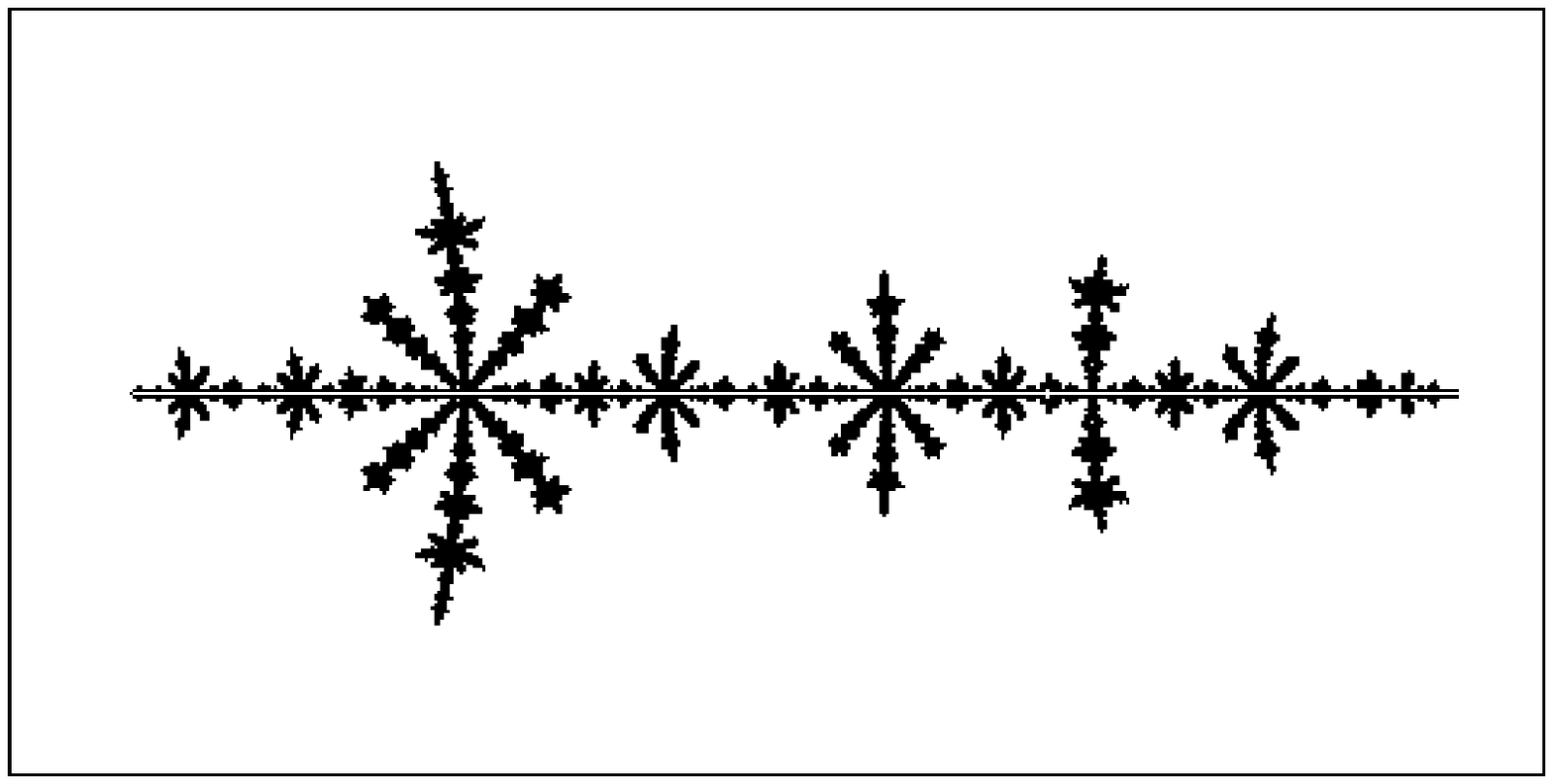 pixels 480 by 240 scaled 550 } 
\smallskip
\centerline{\sl The ``associated" degree 3 polynomial which mimics the critical behaviour of the pair.}  
\bigskip
\vfill\eject

\centerline{\bf References.} 
\bigskip
\bigskip
[DH] A.Douady and J.Hubbard, \'Etude dynamique des polyn\^omes complexes, part I; 
Publ Math. Orsay 1984-1985. 
\smallskip
[M] J.Milnor, Hyperbolic Components in Spaces of Polynomial Maps; Preprint\break 
$\#$1992/3 IMS, SUNY@StonyBrook. 
\smallskip
[P1] A.Poirier, On Postcritically Finite Polynomials; Thesis, SUNY@StonyBrook, 1993 (to appear). 
\smallskip
[P2] A.Poirier, On the realization of Fixed Point Portraits; Preprint $\#$1991/20 IMS, SUNY@StonyBrook. 
\smallskip
[P3] A.Poirier, Hubbard Forest; Preprint $\#$1992/12 IMS, SUNY@StonyBrook. 
\bigskip
\bigskip
 
Typeset in $\TeX$ (July 20, 1994) 
\end